\documentclass[reqno]{amsart}
\usepackage{amssymb,latexsym,amscd}
\theoremstyle{plain}
\newtheorem{theorem}{Theorem}
\newtheorem{lemma}{Lemma}
\newtheorem{proposition}{Proposition}

\theoremstyle{definition}

\newtheorem{remark}{Remark}
\newtheorem{example}{Example}

\numberwithin{equation}{section}

\newcommand{\tP}{\tilde P}
\newcommand{\tX}{\tilde X}

\newcommand{\tG}{\tilde G}

\begin{document}
\begin{abstract}
We study the indefinite metric $G$ in the contact phase space
$(P,\theta)$ of a homogeneous thermodynamical system introduced by
R. Mrugala.  We calculate the curvature tensor, Killing vector
fields, second fundamental form of Legendre submanifolds of $P$ -
constitutive surfaces of different homogeneous thermodynamical
systems. We established an isomorphism of the space $(P,\theta
,G)$ with the Heisenberg Lie group $H_{n}$ endowed with the right
invariant contact structure and the right invariant indefinite
metric. The lift $\tG$ of the metric $G$ to the symplectization
$\tP$ of contact space $(P,\theta)$, its curvature properties, and
its Killing vector fields are studied. Finally we introduce the
"hyperbolic projectivization" of the space $(\tP,{\tilde \theta},
\tG)$ that can be considered as the natural {\bf compactification}
of the thermodynamical phase space $(P,\theta ,G).$
\end{abstract}
\title[On the thermodynamical]%
{The indefinite metric of R.Mrugala and the geometry of the
thermodynamical phase space}
\author{Serge Preston}
\address{Department of Mathematics and Statistics, Portland State
University, Portland, OR, USA} \email{serge@mth.pdx.edu}
\author{James Vargo}
\address{Department of Mathematics, University of Washington, Seattle, WA,USA} \email{jimavargo@hotmail.com}

\maketitle

\footnote{\subjclass[2000]{Primary: 53C50; Secondary:
53D10,74A15}}

\section{Introduction.}
Geometrical methods in the study of homogeneous thermodynamical
systems pioneered by J.Gibbs (\cite{GI}) and C.Caratheodory
(\cite{Car}). They were further developed in the works of
R.Hermann (\cite{H}), R.Mrugala, P.Salamon and their
collaborators, in the dissertations of H.Heemeyer (\cite{He}) and
L.Benayoun (\cite{Be}) to mention just a few. Thermodynamical
metrics (TD-metrics) in the form of the Hessian of a
thermodynamical potential were explicitly introduced by F.Weinhold
(\cite{W}) and, from a different point of view, by G.Ruppeiner
(\cite{R1}).\par

Deeper studies by P.Salamon and his collaborators, by P.Mrugala
and H.Janyszek (see \cite{BSI,J,JM,M1,M2,MNSS,NS,SNI,SAGB,SB},
esp. review papers \cite{DS,M3,R2}) clarified principal properties
of thermodynamical metrics, relations between different
TD-metrics, and their relations to the contact structure of
equilibrium thermodynamical phase space.\par  G. Ruppeiner (see
\cite{R1} and the review \cite{R2}) has developed a covariant
thermodynamical fluctuation theory based on the Riemannian metric
$\eta_{S}$ defined by the second momenta of entropy with respect
to the fluctuations and related the curvature of this metric to
the correlational volume near the critical point.\par

In his work \cite{M2} (see also the review paper \cite{M3})
R.Mrugala introduced the pseudo-Riemannian (indefinite) metric $G$
of signature $(n+1,n)$ in the thermodynamical contact space
$(P,\theta )$ inducing TD metrics on the constitutive surfaces
defined by different thermodynamical potentials (Weinhold metric
for the internal energy and Ruppeiner metric for the entropy).
\par

In the present work we study geometrical properties of the metric
$G$ and its relation to the contact structure $\theta$ of the
thermodynamical phase space (TPS) $(P,\theta)$.\par

In Sec.2 we introduce the (model) thermodynamical phase space $P$
with its contact structure $\theta $. In Section 3 we represent
TPS $P$ as the 1-jet space of the trivial line bundle (Gibbs
space) over the manifold $X$ of {\bf extensive variables} of TPS.
This allows us to define Legendre submanifolds $\Sigma_{\phi}$
(constitutive surfaces of different thermodynamical systems with
given TPS $P$) corresponding to a thermodynamical potential
$\phi\in C^{\infty}(X)$, in their canonical representation (see
\cite{A}).\par

In Section 4 we recall the definitions of the Weinhold and Ruppeiner
metrics.\par

In Section 5 we introduce the indefinite Mrugala metric $G$, a
canonical non-holonomic frame $(X_{i},P_{j},\xi)$.  Then we show the
compatibility of the metric $G$ with the contact structure $\theta$
in a sense that is natural, though different from the conventional
definition used for Riemannian metrics (see \cite{BL}).\par

In Section 6 we determine the Levi-Civita connection $\Gamma$ of
the metric $G$ and work out the formulas (6.4) for covariant
derivatives of the vector fields of the frame $(X_{i},P_{j},\xi).$
\par

In Section 7 we find the Ricci tensor, scalar curvature
$R(G)=\frac{n}{2}$ of metric $G$, and the curvature transformation
$R(X,Y)$ in terms of the frame $(X_{i},P_{j},\xi).$\par

In Section 8 we determine the sectional curvatures of the planes
generated by couples of vectors of the frame
$(\partial_{p_{i}},\partial_{x^{j}},\xi=\partial_{x^{0}}).$ These
curvatures (whenever they are defined) are all zero except for the
plane $P_{i}\wedge
\partial_{x^{i}}$ which has curvature $\frac{3}{4}$.\par

In Section 9 the Lie algebra ${\mathfrak iso}_{G}$ of Killing
vectors of metric $G$ is determined. It is shown that ${\mathfrak
iso}_{G}$ is the Lie algebra ${\mathfrak gl}(n,\mathbb{R})\times
h_{n}$ - the semidirect product of the linear Lie algebra
${\mathfrak gl}(n,\mathbb{R})$ embedded into the symplectic Lie
algebra ${\mathfrak sp}(n,\mathbb{R})$ (with generators
$\{Q^k_l=p_{l}\partial_{p_{l}}-x^k
\partial_{x^{k}}\}$) and of the Heisenberg Lie algebra ${\mathfrak
h}_{n}$ with generators $\{\xi
=\partial_{x^{0}},A_i=\partial_{p_{i}}+x^{i}\partial_{x^{0}},B_j=-\partial_{x^{j}}\}$
and commutative relations (9.3-4).\par

In Section 10 we calculate the second fundamental form of
constitutive surfaces (Legendre submanifolds) of contact manifold
$(P,\theta)$.\par

In Sec.11 the constitutive hypersurface
$x^0+\sum_{i}p_{i}x^{i}=0$, defined by the homogeneity condition
of thermodynamical potentials (see \cite{CA}) is introduced and
studied.\par

In Section 12 we establish the isomorphism of the TPS $(P,\theta
,G)$ with the Heisenberg group $H_{n}$ endowed with the right
invariant contact structure and right invariant indefinite metric.
More specifically, we prove the following
\begin{theorem}
The diffeomorphism $\chi $ defined by
\[
\chi :\ g=\begin{pmatrix} 1 & {\bar a} & z\\ 0 & I_{n} & {\bar b}\\
0 & 0 & 1
\end{pmatrix} \rightarrow m=\begin{pmatrix} -z\\{\bar b}\\{\bar
a}\end{pmatrix}.
\]
determines an isomorphism of
the "thermodynamical metric contact manifold" $(P,\theta ,G)$ with
the Heisenberg group $H_{n}$ endowed with the right invariant
contact from $\theta_{H}$ and the right invariant metric $G_{H}$ of
signature $(n+1,n).$
\end{theorem}

In Part II of the paper we study the {\bf symplectization} $(\tP,
{\tilde \theta}, \tG)$ of the contact manifold $(P,\theta ,G)$ with
the indefinite metric $G$. We lift $G$ to the indefinite (of
signature (n+1,n+1)) metric $\tG$, calculate the Levi-Civita
connection, the Ricci tensor, and the scalar curvature
$R(\tG)=(n+1)(n+2)$ of $\tG$.  We get, in particular, that
$Ric(\tG)=\frac{n+2}{2}\tG$, i.e. $\tG$ is an {indefinite
(pseudo-Riemannian) Einstein metric}.\par

In Section 16 and in the Appendix we determine the Lie algebra
${\mathfrak iso}_{\tG}$ of the Killing vector fields of metric
$\tG$. We prove that ${\mathfrak iso}_{\tG}\simeq {\mathfrak
sl}(n+2,\mathbb{R})$.\par

In Section 19 we show that the contact metric space $(P,\phi,
\theta ,\xi)$ endowed with the almost contact structure given by
the (1,1)-tensor (see \cite{M2})
\[
\phi =\begin{pmatrix} 0 & p_{1},\ldots ,p_{n} & 0\\
0 & 0_{n} & I_{n}\\
0 & -I_{n} & 0_{n}
\end{pmatrix}
\]
is a {\bf Sasakian metric} in that the almost complex structure $J$
defined on the space $\tP$ by the formula
\[
J(X,f\partial_{p_{0}})=(\phi (X)-f\xi ,\eta(X))\partial_{p_{0}})
\]
for any function $f\in C^{\infty}(\tP)$ {\bf is integrable}.\par

In Section 20 we construct the symplectomorphism of the "positive
quadrant"
\[
(\tP_{+}=\{(p_{i},x_{i})\in \tP \vert p_{i}>0\},d{\tilde \theta })
\]
of the manifold $(\tP, {\tilde \theta }$ is isomorphic to the
product $\prod_{i=0}^{i=n}A_{1}^{i}$ of (n+1) copies of the affine
group
\[
A_{1}=\{ \begin{pmatrix} e^t & c\\ 0 & 1 \end{pmatrix} \vert
t,c\in \mathbb{R} \}
\]
 endowed with the symplectic structure
generated by the right invariant 1-form
$\sum_{i=0}^{i=n}(dz_{i}-z_{i}h_{i}^{-1}dh_{i})$ on the Lie group
$\prod_{i=0}^{i=n}A_{1}^{i}$. \par

Finally, in Section 21 we define the "hyperbolic projectivization"
$({\hat P}\simeq P_{2n+1}(\mathbb{R}) ,{\hat \theta},{\hat G})$ of
the space $(\tP ,{\tilde \theta},\tG)$ - a natural
compactification of the TPS space $(P,\theta ,G).$ This
projectivization will be used for the study of geometrical
properties of thermodynamical systems in the continuation of this
work.\par
\vskip1cm

\centerline{ \textbf{Part I.}}

\section{The contact structure of homogeneous thermodynamics.}
A phase space of the Homogeneous Thermodynamics (thermodynamical
phase space, or TPS) is the (2n+1)-dimensional vector space
$P=\mathbb{R}^{2n+1}$ endowed with the standard contact structure
(\cite{AG,H})

\begin{equation} \theta
=dx^{0}+\sum_{i=1}^{n}p_{l}dx^{l}.
\end{equation}

The horizontal distribution of this structure is generated by two
families of vector fields: $P_{i},X_{i}$

\[
D_{m}=<P_{i}=\partial_{p_{l}},\
X_{i}=\partial_{x^{i}}-p_{i}\partial_{x^{0}}>.
\]
The 2-form
\[
\omega =d\theta =\sum_{i=1}^{n}dp_{l}\wedge dx^{l}
\]
is a nondegenerate, symplectic form on the distribution $D$.\par

The Reeb vector field, uniquely defined as the generator $\xi $ of
$ker(d\theta)$ satisfying $\theta (\xi )=1$, is
simply
\[
\xi = \partial_{x^{0}}.
\]

\section{Gibbs space. Legendre surfaces of equilibrium.}
Constitutive surfaces of concrete thermodynamical systems are
determined by their "constitutive equations", which, in their
fundamental form determine the value of a {\bf thermodynamical
potential} $x^{0}=E(x^{i})$ as the function of n extensive variables
$x^{i}$. Dual intensive variables are determined then as the partial
derivatives of the thermodynamical potential by the extensive
variables: $p_{i}=\frac{\partial E}{\partial x^{i}}$.\par
  Thus, a constitutive surface represents
the {\bf Legendre submanifold} (maximal integral submanifold)
$\Sigma_{E}$ of the contact form $\theta$ projecting
diffeomorphically to the space $X$ of variables $x^{i}$. Space $Y$
of variables $x^{0},x^{i},\ i=1,\ldots ,n$ is, sometimes, named the
{\bf Gibbs space} of the thermodynamical potential $E(x^{i})$.
Thermodynamical phase space $(P,\theta )$ (or, more precisely, its
open subset) appears as the {\bf first jet space}
$J^{1}(Y\rightarrow X)$ of the (trivial) line bundle
$\pi:Y\rightarrow X$. Projection of $\Sigma_{E}$ to the Gibbs space
$Y$ is the {\bf graph} $\Gamma_{E}$ of the constitutive law
$E=E(x^{i})$.
\par

Another choice of the thermodynamical potential together with the
n-tuple of extensive variables leads to another representation of
an open subset of TPS $P$ as the 1-jet bundle of the corresponding
Gibbs space.  It is known that any Legendre submanifold of a
contact form $\theta $ can locally be presented in this form for
some choice of the set of extensive variables and of
thermodynamical potential as the function of these variables.\par

We will be using a local description of Legendre submanifolds
which takes a slightly different form. Let $P^{2n+1}$ be a contact
manifold. The following result characterizes (locally) all
Legendre submanifolds (V.I. Arnold, \cite{A,AG}).\par

Choose (local) Darboux coordinates $(x^{0},x^{i},p_{j})$ in which
$\theta= dx^0+p_k dx^k$. Let $I,J$ be a partition of the set of
indices $1,\dots,n$, and consider any function $\phi(p_i,x^j),\ i\in
I,j\in J$. Then the following equations define a Legendre
submanifold $\Sigma_{\phi}$:
\begin{equation}
\begin{cases}
x^0&=\phi -\sum_{i\in I} p_i\frac{\partial \phi}{\partial
p_{i}},\\
p_j&=-\frac{\partial \phi}{\partial x^{j}},\ j\in J,\\
x^i&=\frac{\partial \phi}{\partial p_i},\ i\in I.
\end{cases}
\end{equation}

\medskip

Moreover, every Legendre submanifold is locally given by some choice
of a splitting $I,J$ and of a function $\phi(p_i,i\in I,x^j,j\in J)
$.\par

In physics, the most commonly used thermodynamical potentials are:
internal energy, entropy, free energy of Helmholtz, enthalpy and
the free Gibbs energy.\par

On the intersection of the domains of these representations,
corresponding points are related by a {\bf Legendre transformations}
(see \cite{A}).\par
\begin{example}
As an example of such a thermodynamical system, consider the van
der Waals gas - a system with two thermodynamical degrees of
freedom. Space $P$ is 5-dimensional (for 1 mole of gas) with the
canonical variables $(U,(T,S),(-p,V))$ (internal
energy,temperature, entropy, -pressure, volume), the contact form
\[
\theta =dU-TdS+pdV,
\]
and the fundamental constitutive law
\[
U(S,V)=(V-b)^{\frac{R}{C_{V}}}e^{\frac{S}{c_{V}}}-\frac{a}{V},
\]
where $R$ is the Ridberg constant, $c_{V}$ is the heat capacity at
constant volume, $a,b$ are parameters of the gas reflecting the
interaction between molecules and the part of volume occupied by
molecules respectively.
\end{example}

\section{Thermodynamical metrics of Weinhold and Ruppeiner.}
A thermodynamical metric $g_{U}$ (the Weinhold metric) in the space
$X$ of extensive variables corresponding to the choice of internal
energy $U$ as the thermodynamical potential $E$ was explicitly
introduced by F.Weinhold (see \cite{W}) as the Hessian of the
internal energy $U(x^{i})$
\begin{equation}
g_{U\ ij}=\frac{\partial^{2}U}{\partial x^{i}\partial x^{j}},
\end{equation}
\par
G. Ruppeiner's metric $g_{S}$ corresponding to the choice of entropy
$S$ as the thermodynamical potential $E$ was defined by the same
formula and intensively studied by Ruppeiner (\cite{R1}) in the
framework of the fluctuational theory of thermodynamical
systems.\par

Interest in these metrics is partly due to the fact that the
definiteness of $g_{U}$ (positive or negative) at a point $x\in X$
of the of the constitutive surface delivers the local criteria of
stability of the equilibria given by the corresponding point of the
surface $\Gamma_{E}$ or $\Sigma_{E}$ (see \cite{CA}). \par

Later both of these metrics were studied by P.Salomon and his
collaborators and by R.Mrugala and H.Janyszek (see ref.).
Geometrical properties of these metrics were studied for TD
systems with small degree of freedom (small n). An interesting and
important in applications meaning was assigned to the length of
curves ("processes") in the space $X$ (\cite{DS,SAGB,SB}).  It was
suggested that the curvature of these metrics is related to the
interactions in the (microscopical) system, and that singularities
of the scalar curvature of these metrics were related to the
properties of system near the phase transition and the triple
point of the thermodynamical system (\cite{R1,R2}).

\section{Indefinite thermodynamical metric $G$ of R. Mrugala.}
In the paper \cite{M2}, R.Mrugala defined a pseudo-Riemannian
(indefinite) metric $G$ of signature $(n+1,n)$ in the
thermodynamical contact space $(P,\theta ).$  It is given by the
formula:  \begin{equation} G=2dp_k \odot dx^k+\theta\otimes\theta,
\end{equation}
where $dp_k
\odot dx^k=\frac{1}{2}(dp_k\otimes dx^k+dx^k\otimes dp_k)$ is the
symmetrical product of 1-forms.
\par\bigskip

Its physical motivation is two-fold.  First, it was derived by means
of statistical mechanics.  Second, its reduction to the Legendre
submanifolds $\Sigma_{\phi}$ corresponding to the choice of entropy
or internal energy as the TD potential coincides with the previously
studied Ruppeiner and Weinhold metrics $\eta_{\phi}$.  Indeed, let
$\phi=S$ or $U$, and form the Legendre submanifold $\Sigma_\phi$ of
the following type:

\[x^0=\phi(x^1,\dots,x^n) \hspace{.25 in}
p_k=\frac{\partial\phi}{\partial x^k}.\]

In the coordinates $\{x^k\}$, the restriction of $G$ to
$\Sigma_\phi$ has components
\[
g_{ij}=\frac{\partial^2\phi}{\partial x^i
\partial x^j}.
\]

In coordinates $x^{0}; p_{s};x^{i}$, metric $G$ is given by the
following matrix:
\begin{equation}
G=(G_{ij})=
\begin{pmatrix}
1 &  0  &  p_{j}\\
0 &  0  &  I_{n}\\
p_{i} & I_{n} & p_{i}p_{j}
\end{pmatrix} .
\end{equation}
It is easy to see that
\begin{equation}
det(G)=(-1)^{n},
\end{equation}
therefore metric $G$ is non-degenerate.  \par
\par
The inverse (covariant) metric to $G$ is given by
\begin{equation}
G^{-1}=(G^{ij})=
\begin{pmatrix}
1 &  -p_{j}  &  0\\
-p_{i} &  0  &  I_{n}\\
0 & I_{n} & 0
\end{pmatrix} .
\end{equation}

\subsection{Non-holonomic frame $(\xi, P_{i},X_{i})$.}

It is convenient to introduce the following non-holonomic frame of
the tangent bundle $T(P)$ (see \cite{Ki})
\begin{equation}
\xi=\partial_{x^{0}},\ P_{l}=\partial_{p_{l}},\
X_{i}=\partial_{x^{i}}-p_{i}\partial_{x^{0}}
\end{equation}
whose only nonzero commutator relation is

\[
[P_{i},X_{j}]=-\delta_{ij}\xi .
\]
In this frame the metric $G$ takes the simple form
\begin{equation}
{G}=(G_{ij})=
\begin{pmatrix}
1 &  0  &  0\\
0 &  0  &  I_{n}\\
0 & I_{n} & 0
\end{pmatrix} .
\end{equation}
The dual coframe of the frame (5.5) is given by
\begin{equation}
(\theta , dp_{l},dx^{i}).
\end{equation}

\par

Using the frame (5.5) one can easily find positive and negative
sub-bundles of tangent bundle $T(P)$ - they are generated,
correspondingly, by the tangent vectors
\begin{equation}
T_{+}=<\partial^{x^{0}},
\frac{1}{2}(\partial_{p_{l}}+\partial_{x^{l}}-p_{l}\partial_{x^{0}})>,\
T_{-}=
<\frac{1}{2}(\partial_{p_{l}}-\partial_{x^{l}}+p_{l}\partial_{x^{0}})>.
\end{equation}
In this basis, metric $G$ takes the standard indefinite form
\begin{equation}
{\bar G}=(G_{ij})=
\begin{pmatrix}
1 &  0  &  0\\
0 &  I_{n}  &  0\\
0 & 0 & -I_{n}
\end{pmatrix} .
\end{equation}
In terms of this frame the light cone at each point $m\in P$ is
given by the standard quadric, i.e. $X=f\xi
+g^{k}P_{k}+h^{s}X_{s}$ belongs to the light cone at a point $m$
if and only if
\[
f^{2}+2\sum_{k}g^{k}h^{k}=0.
\]

\subsection{Compatibility of contact structure and metric $G$.}
The notion of compatibility between an almost contact structure and
a Riemannian metric that was introduced by Sasaki has become a
classical notion (see \cite{BL}). We recall here that an {\bf almost
contact structure} on a manifold $M^{2n+1}$ is defined by a triple
$(\theta, \xi, \phi )$ of a 1-form $\theta$, a vector field $\xi $
and a (1,1)-tensor field $\phi $ satisfying the conditions
\begin{equation}
\phi^{2}=-I+\theta\otimes \xi ;\ \theta(\xi )=1.
\end{equation}
From these properties the subsequent relations follow
\[
\phi (\xi )=0;\ \theta \circ \phi =0;\ rank(\phi )=2n.
\]
\par
A Riemannian metric $g$ is said to be {\bf compatible with the
almost contact structure $(\theta ,\xi , \phi )$} (or associated
with it) if
\begin{equation}
g(\phi X,\phi Y)=g(X,Y)-\theta(X)\theta(Y)
\end{equation}
for all tangent vectors $X,Y$ at all points $m\in M$ (see
\cite{BL}).\par

It is known that any almost contact structure admits a (far from
being unique) compatible Riemannian metric,\cite{BL}.\par
For a contact manifold $(M\theta )$ the condition of
compatibility above is equivalent to the following two conditions
taken together: 1) on the contact distribution $D=ker(\theta)$
$\phi$ is a $g$-orthogonal transformation, and 2) $\xi $ is
$g$-orthogonal to $D$.\par

In a case of an {\bf indefinite metric} $g$ on
$(M^{2n+1},\theta)$, the compatibility condition (5.11) should be
modified if we would like to incorporate even the most simple
indefinite metric
\[
g=\begin{pmatrix} 1& 0 & 0\\ 0 & 0 & 1 \\0 & 1 & 0\end{pmatrix}
\]
defined in standard 3D contact space $(\mathbb{R}^3,
x^{0},p_{1},x^{1})$ with the contact form $\theta
=dx^{0}+p_{1}dx^{1}$.  It is sufficient to check the condition for
two basic horizontal vectors
$X=\partial_{x^{1}}-p_{1}\partial_{x^{0}},\ P=\partial_{p_{1}}.$
\par

We have $g(X,Y)=1$ and $\phi(X)=-X,\phi(P)=P$, so that $g(\phi
X,\phi Y)=-1.$  Thus, condition (5.11) is not fulfilled.\par

For the Mrugala metric, the $(1,1)$-tensor $\phi$ of the
associated almost contact structure has the form (see \cite{M2})
\begin{equation}
\phi =\begin{pmatrix} 0 & p_{1},\ldots ,p_{n} & 0\\
0 & 0_{n} & I_{n}\\
0 & -I_{n} & 0_{n}
\end{pmatrix}.
\end{equation}

We have
\[
\phi (\xi)=0, \phi (X_{i})=P_{i}, \phi (P_{k})=-X_{k}
\]
and
\[
G(\phi(X_{i}),\phi(P_{j}))=G(P_{i},-X_{j})=-\delta_{ij},
\]
while $G(X_{i},P_{j})=\delta_{ij}$.  Thus, the Mrugala metric does
not satisfy the conventional compatibility condition with the
contact structure.\par

In the forthcoming paper \cite{P} we analyze the situation of an
indefinite metric determined on a contact manifold and suggest a
definition of compatibility suited for indefinite metrics of
arbitrary signature.  In particular, for the standard 3D metric
above and for the Mrugala metric $G$ this definition reads as
\begin{equation}
g(\phi X,\phi Y)=-g(X,Y)-\theta(X)\theta(Y),
\end{equation}
which is satisfied by both of these metrics.

\section{Levi-Civita Connection.}
In this section, we will compute the Christoffel connection
coefficients $\Gamma^{\alpha}_{\beta\gamma}$ and the covariant
derivatives of vector fields of the frame (5.5) with respect to
the vector fields of the same frame.\par

  For calculation of the Christoffel coefficients we define combinations
\[\{\alpha\beta,\,\gamma\}=G_{\alpha\gamma,\,\beta}+G_{\beta\gamma,\,\alpha}-G_{\alpha\beta,\,\gamma}\]
symmetrical by $(\alpha \beta)$.\par With this notation, the
formula for the connection coefficients is (\cite{KN})
\[\Gamma^{\alpha}_{\beta\gamma}=\frac{1}{2}G^{\alpha
s}\{\beta\gamma,\, s\}\]

First, note that the only nonzero derivatives of the metric G are
given by:  \[G_{x^{i}0, p_{j}}=G_{0x^{i},\, p_{j}}=\delta_{ij},\
G_{x^{i}x^{j}, p_{k}}=\delta_{ki}p_{j}+\delta_{kj}p_{i}.\]

This implies that the only nonzero combinations
$\{\alpha\beta,\,\gamma\}$ are the following ones:
\[\{0x^{i},\,p_{j}\}=-\delta_{ij};\ \{0p_{j},\,x^{i}\}=\delta_{ij};\ \{x^{i}p_{j},\,0\}=\delta_{ij};\]
\[\{x^{i}x^{j},\,p_{k}\}=-\left(\delta_{ik}p_{j}+\delta_{jk}p_{i}\right);
\ \{x^{i}p_{j},\,x^{k}\}=\delta_{ij}p_{k}+\delta_{kj}p_{i}.\]

Now we will calculate the connection coefficients.\par

For any $\alpha$, $\Gamma^{\alpha}_{00}=0$.\par

Next, $\Gamma^{0}_{0\gamma}=-\frac{1}{2}p_{k}\{0\gamma,\,p_{k}\}$.
And this quantity is 0 unless $\gamma=x^{i}$, in which case, we
obtain $\Gamma^{0}_{0x^{i}}=\frac{1}{2}p_{i}$.\par

In the next case, we consider $\beta\neq 0$ and $\gamma\neq
0$\[\Gamma^{0}_{\beta\gamma}=\frac{1}{2}{[\{\beta\gamma,
0\}-p_{k}\{\beta\gamma, p_{k}\}]}.\]
\par
So we obtain:  \[\Gamma^{0}_{x^{i}p_{j}}=\frac{1}{2}\delta_{ij},\]
\[\Gamma^{0}_{x^{i}x^{j}}=\frac{1}{2}p_{k}\left(\delta_{ik}p_{j}+\delta_{jk}p{i}\right)=p_{i}p_{j},\]
\[\Gamma^{0}_{p_{i}p_{j}}=0.\]

The next case to consider is when $\beta=0$, and the other indices
are nonzero. Then we have:

\[\Gamma^{\alpha}_{0\gamma}=\frac{1}{2}{[G^{\alpha p_{k}}\{0\gamma,
p_{k}\}+G^{\alpha x^{l}}\{0\gamma, x^{l}\}}].\]

If $\alpha=p_{i}$, then we get
$\frac{1}{2}\delta^{i}_{l}\{0\gamma, x^{l}\}=\frac{1}{2}\{0\gamma,
x^{i}\}$.  This is $\frac{1}{2}\delta^{i}_{j}$ if
$\gamma=p_{j}$,and is zero otherwise.

If $\alpha=x^{i}$, then we get
$\frac{1}{2}\delta^{i}_{k}\{0\gamma, p_{k}\}=\frac{1}{2}\{0\gamma,
p_{i}\}.$ If $\gamma=x^{j}$, then this is equal to
$-\frac{1}{2}\delta^{i}_{j}.$  Otherwise it is zero.

The next case is when the upper index is $p_{i}$, and the lower
indices are nonzero.
\[\Gamma^{p_{i}}_{\beta\gamma}=\frac{1}{2}{[-p_{i}\{\beta\gamma,
0\}+\{\beta\gamma, x^{i}\}]}.\] The only nonzero term comes when
$\beta\gamma=x^{j}p_{k}.$  So we obtain:
\[\Gamma^{p_{i}}_{x^{j}p_{k}}=\frac{1}{2}{[-p_{i}\delta^{j}_{k}+\delta^{j}_{k}p_{i}+\delta^{i}_{k}p_{j}]}=\frac{1}{2}\delta^{i}_{k}p_{j}.\]

The final case is the one in which the upper index is $x^{i}$ and
the lower indices are nonzero.
\[\Gamma^{x^{i}}_{\beta\gamma}=\frac{1}{2}\{\beta\gamma,
p_{i}\}\]This is only nonzero if $\beta\gamma=x^{j}x^{k}$, and in
this case, we get
\[\Gamma^{x^{i}}_{x^{j}x^{k}}=-\frac{1}{2}\left(\delta^{i}_{j}p_{k}+\delta^{i}_{k}p_{j}\right).\]

To summarize, the {\bf nonzero} Christoffel coefficients are given
by:
\begin{multline}
\Gamma^{0}_{0x^{i}}=\frac{1}{2}p_{i};\
\Gamma^{0}_{x^{i}p_{j}}=\frac{1}{2}\delta_{ij};\
\Gamma^{0}_{x^{i}x^{j}}=p_{i}p_{j};\
\Gamma^{p_{i}}_{0p_{j}}=\frac{1}{2}\delta^{i}_{j};\
\Gamma^{x^{i}}_{0x^{j}}=-\frac{1}{2}\delta^{i}_{j};\\
\Gamma^{p_{i}}_{x^{j}p_{k}}=\frac{1}{2}\delta^{i}_{k}p_{j};\
\Gamma^{x^{i}}_{x^{j}x^{k}}=-\frac{1}{2}\left(\delta^{i}_{j}p_{k}+\delta^{i}_{k}p_{j}\right).
\end{multline}

Finally, we calculate the trace 1-form $\gamma$ of the connection
$\Gamma$ (see \cite{Po}), whose components in the coordinate coframe
$dx^{j}, dp_{j}$ are given by
$\gamma_{\alpha}=\Gamma^{\beta}_{\alpha\beta}.$
\[\gamma_{0}=\Gamma^{p_{i}}_{0p_{i}}+\Gamma^{x^{i}}_{0x^{i}}=\frac{1}{2}\delta^{i}_{i}-\frac{1}{2}\delta^{i}_{i}=0\]
\[\gamma_{p_{j}}=\Gamma^{0}_{p_{j}o}+\Gamma^{x^{i}}_{p_{j}x^{i}}+\Gamma^{p_{l}}_{p_{j}p_{l}}=0\]
\[\gamma_{x^{j}}=\Gamma^{0}_{x^{j}0}+\Gamma^{p_{i}}_{x^{j}p_{i}}+\Gamma^{x^{i}}_{x^{j}x^{i}}=\frac{1}{2}p_{j}+\frac{1}{2}\delta^{i}_{i}p_{j}-\frac{1}{2}\left(\delta^{i}_{j}p_{i}+\delta^{i}_{i}p_{j}\right)=0\]

So  in the coframe $\theta , dx^{j}, dp_{j}$, $\gamma=0$ (recall
that the form $\gamma $ depends on the choice of a coordinate system
- in other coordinates, the form $\gamma $ differs from zero by the
differential $df$ of a function (see \cite{Po}).
\par
\vskip0.5cm

It is helpful for the future calculations to
point out that the covariant derivative of the connection $\Gamma$
takes a particularly simple form when expressed in the
non-holonomic basis:
\[\xi=\partial_{x^0}, P_i=\partial_{p_i},
X_i=\partial_{x^i}-p_i\xi,\] defined above. \par

Let $G$ be any pseudo-Riemannian metric.  Then, for any vector
fields $X,Y,Z$, the Levi-Civita connection of metric $G$
satisfies:
\[YG(X,Z)=G(\nabla_{Y}X,Z)+G(X,\nabla_{Y}Z);\]
\[ZG(X,Y)=G(\nabla_{Z}X,Y)+G(X,\nabla_{Z}Y);\]
\[XG(Y,Z)=G(\nabla_{X}Y,Z)+G(Y,\nabla_{X}Z).\]
If we add the first two equations and subtract the third, the
result is:
\begin{equation}
2G(X,\nabla_{Y}Z)=YG(X,Z)+ZG(X,Y)-XG(Y,Z)+G(Y,[X,Z])+G(Z,[X,Y])+G(X,[Y,Z]).
\end{equation}

Here we have used the fact that the Levi-Civita connection is
symmetric; that is, $\nabla_XY-\nabla_YX=[X,Y]$ for all
$X,Y$.\bigskip

In taking $X,Y,Z$ from  the vectors of the frame
$\{\xi,P_i,X_j\}$, scalar products $G(X,Y),G(Y,Z),G(Z,X)$ will all
be constant; therefore, the first three terms on the right side of
(6.2) will vanish, leaving:

\begin{equation} 2G(X,\nabla_{Y}Z)=G(Y,[X,Z])+G(Z,[X,Y])+G(X,[Y,Z]).
\end{equation}

Among the basic vectors, the only pair of vectors with non-zero
Lie-bracket is $[P_i,X_j]=-\delta_{ij}\xi$. It follows that if we
substitute basic vectors into the equation above, the right side
will equal zero unless two of the vectors are $P_i$ and $X_j$,
respectively. Since their bracket is proportional to $\xi$, which is
orthogonal to the contact distribution, the third vector must be
$\xi$.  In particular, we immediately obtain the following
relations:
\[\nabla_\xi\xi=0, \nabla_{P_i}P_j=0, \nabla_{X_i}X_j=0.\]

Additionally, we see that the only nonzero component of
$\nabla_{P_i}X_j$ is the $\xi$ component, which is found by:
\[2G(\xi,\nabla_{P_i}X_j)=G(\xi,\delta_{ij}\xi)=-\delta_{ij}.\]

Therefore, $\nabla_{P_i}X_j=-\frac{1}{2}\delta_{ij}\xi$.  Note
that interchanging $P_i$ with $X_j$ changes the sign in the right
side: $\nabla_{X_j}P_i=\frac{1}{2}\delta_{ij}\xi$.\bigskip

The next equation to consider is:
\[2G(P_i,\nabla_{\xi}X_j)=G(\xi,[P_i,X_j])=-\delta_{ij}.\]
It follows that $\nabla_{\xi}X_j=-\frac{1}{2}X_j$.  Interchanging
the roles of $P_i,X_j$ yields $\nabla_{\xi}P_i=\frac{1}{2}P_i$. By
the symmetry of the connection, $\nabla_{X_j}\xi=\nabla_{\xi}X_j$
and $\nabla_{P_i}\xi=\nabla_{\xi}P_i$ since $[\xi ,P_{i}]=0$.
\bigskip

In summary, the curvature derivatives of the connection $\Gamma $
in the (canonical) frame (5.5) are given by the following
equations:
\begin{equation}
\begin{aligned}
\nabla_\xi\xi &=0, \nabla_{P_i}P_j =0, \nabla_{X_i}X_j
=0;\\
\nabla_{\xi}P_i&=\nabla_{P_i}\xi=\frac{1}{2}P_i;\\
\nabla_{\xi}X_j&=\nabla_{X_j}\xi=-\frac{1}{2}X_j;\\
-\nabla_{P_i}X_j&=\frac{1}{2}\delta_{ij}\xi=\nabla_{X_i}P_j.
\end{aligned}
\end{equation}

\section{Ricci and scalar curvatures.}
In this section, we calculate the Ricci tensor and the full
curvature tensor in the form of a transformation of the tangent
bundle.

Recall the formula for the components of the Ricci tensor
(\cite{KN}):
\begin{equation}
R_{\alpha\beta}=\Gamma^{\mu}_{\beta\alpha,
\mu}-\Gamma^{\mu}_{\mu\alpha,
\beta}+\Gamma^{\mu}_{\mu\gamma}\Gamma^{\gamma}_{\beta\alpha}-\Gamma^{\mu}_{\beta\gamma}\Gamma^{\gamma}_{\mu\alpha}
\end{equation}

The Ricci tensor of a Riemannian space is symmetric.  In addition to
this, by the remarks above, the result of contracting an upper and
lower index of the Christoffel coefficients $\Gamma^{i}_{jk}$ is
zero for our metric $G$. Therefore the middle two terms in the right
side of the above formula are zero. So for the present calculations,
we may use the formula:

\[R_{\alpha\beta}=\Gamma^{\mu}_{\beta\alpha,
\mu}-\Gamma^{\mu}_{\beta\gamma}\Gamma^{\gamma}_{\mu\alpha};\]

\[R_{00}=-\Gamma^{\mu}_{0\gamma}\Gamma^{\gamma}_{\mu0}=-\left(\Gamma^{\mu}_{0p_{i}}\Gamma^{p_{i}}_{\mu0}+\Gamma^{\mu}_{0x^{i}}\Gamma^{x^{i}}_{\mu0}\right)
=-\left(\frac{1}{4}\delta^{j}_{i}\delta^{i}_{j}+\frac{1}{4}\delta^{j}_{i}\delta^{i}_{j}\right)=-\frac{n}{2};\]

\[R_{0p_{i}}=-\Gamma^{\mu}_{p_{i}\gamma}\Gamma^{\gamma}_{\mu0}=
-\left(\Gamma^{\mu}_{p_{i}0}\Gamma^{0}_{\mu0}+\Gamma^{\mu}_{p_{i}x^{j}}\Gamma^{x^{j}}_{\mu0}\right)=0;\]

\[R_{0x^{i}}=-\Gamma^{\mu}_{x^{i}\gamma}\Gamma^{\gamma}_{\mu0}=-\left(\Gamma^{p_{j}}_{x^{i}\gamma}\Gamma^{\gamma}_{p_{j}0}+\Gamma^{x^{j}}_{x^{i}\gamma}\Gamma^{\gamma}_{x^{j}0}\right)=-\left(\frac{1}{2}\delta^{j}_{k}p_{i}\frac{1}{2}\delta^{k}_{j}+\Gamma^{x^{j}}_{x^{i}\gamma}\Gamma^{\gamma}_{x^{j}0}\right)=\]
\[=-\left(\frac{1}{4}np_{i}+\Gamma^{x^{j}}_{x^{i}x^{k}}\Gamma^{x^{k}}_{x^{j}0}+\Gamma^{x^{j}}_{x^{i}0}\Gamma^{0}_{x^{j}0}\right)=-\left(\frac{n}{4}p_{i}+\frac{1}{4}\delta^{k}_{j}\left(\delta^{j}_{i}p_{k}+\delta^{j}_{k}p_{i}\right)+-\frac{1}{4}\delta^{j}_{i}p_{j}\right)=\]
\[=-\left(\frac{n}{4}p_{i}+\frac{1}{4}{[p_{i}+np_{i}]}-\frac{1}{4}p_{i}\right)=-\frac{n}{2}p_{i};\]

\[R_{p_{i}p_{j}}=-\Gamma^{\mu}_{p_{j}\gamma}\Gamma^{\gamma}_{\mu
p_{i}}=-\left(\Gamma^{\mu}_{p_{j}0}\Gamma^{0}_{\mu
p_{i}}+\Gamma^{\mu}_{p_{j}x^{k}}\Gamma^{x^{k}}_{\mu
p_{i}}\right)=-\Gamma^{0}_{p_{j}x^{k}}\Gamma^{x^{k}}_{0p_{i}}=0;\]

\[R_{p_{i}x^{j}}=\frac{\partial}{\partial
p_{k}}\left(\frac{1}{2}\delta^{k}_{i}p_{j}\right)-\Gamma^{\mu}_{x^{j}\gamma}\Gamma^{\gamma}_{\mu
p_{i}}=\frac{1}{2}\delta_{ij}-\left(\Gamma^{0}_{x^{j}\gamma}\Gamma^{\gamma}_{0p_{i}}+\Gamma^{x^{k}}_{x^{j}\gamma}\Gamma^{\gamma}_{x^{k}
p_{i}}\right)=\]
\[=\frac{1}{2}\delta_{ij}-\left(\frac{1}{4}\delta_{jk}\delta^{k}_{i}
+\Gamma^{x^{k}}_{x^{j}0}\Gamma^{0}_{x^{k}p_{i}}\right)=\frac{1}{2}\delta_{ij}-\left(\frac{1}{4}\delta_{ij}
-\frac{1}{4}\delta^{k}_{j}\delta_{ki}\right)=\frac{1}{2}\delta_{ij};\]

\[R_{x^{i}x^{j}}=-\Gamma^{\mu}_{x^{j}\gamma}\Gamma^{\gamma}_{\mu
x^{i}}=-\left(\Gamma^{\mu}_{x^{j}0}\Gamma^{0}_{\mu
x^{i}}+\Gamma^{\mu}_{x^{j}p_{k}}\Gamma^{p_{k}}_{\mu
x^{i}}+\Gamma^{\mu}_{x^{j}x^{k}}\Gamma^{x^{k}}_{\mu
x^{i}}\right)=\]
\[=-\left(\Gamma^{0}_{x^{j}0}\Gamma^{0}_{0x^{i}}+\Gamma^{x^{k}}_{x^{j}0}\Gamma^{0}_{x^{k}x^{i}}+\Gamma^{p_{r}}_{x^{j}p_{k}}\Gamma^{p_{k}}_{p_{r}x^{i}}+\Gamma^{0}_{x^{j}x^{k}}\Gamma^{x^{k}}_{0x^{i}}+\Gamma^{x^{r}}_{x^{j}x^{k}}\Gamma^{x^{k}}_{x^{r}x^{i}}\right)=\]
\[=-\left(\frac{1}{4}p_{i}p_{j}-\frac{1}{2}\delta^{k}_{j}p_{k}p_{i}+\frac{1}{4}\delta^{r}_{k}p_{j}\delta^{k}_{r}p_{i}+-\frac{1}{2}\delta^{k}_{i}p_{j}p_{k}+\frac{1}{4}\left(\delta^{r}_{j}p_{k}+\delta^{r}_{k}p_{j}\right)\left(\delta^{k}_{r}p_{i}+\delta^{k}_{i}p_{r}\right)\right)=\]
\[=-\left(-\frac{3}{4}p_{i}p_{j}+\frac{n}{4}p_{i}p_{j}
+\frac{1}{4}\left(p_{i}p_{j}+p_{i}p_{j}+np_{i}p_{j}+p_{i}p_{j}\right)\right)=-\frac{n}{2}p_{i}p_{j}.\]

As a result, the components of the symmetric Ricci tensor are given
by:
\[R_{00}=-\frac{n}{2},\
R_{0p_{i}}=0,\ R_{0x^{i}}=-\frac{n}{2}p_{i},\]
\[R_{p_{i}p_{j}}=0,\
R_{p_{i}x^{j}}=\frac{1}{2}\delta_{ij},\
R_{x^{i}x^{j}}=-\frac{n}{2}p_{i}p_{j}.\]
\par

To calculate the scalar curvature $G^{\alpha\beta}R_{\beta\alpha}.
$, we contract the Ricci tensor

\[\mathcal{R}=G^{\alpha0}R_{0\alpha}+G^{\alpha
p_{i}}R_{p_{i}\alpha}+G^{\alpha x^{i}}R_{x^{i}\alpha}=\]
\[=\left(R_{00}-p_{k}R_{0p_{k}}\right)+\left(-p_{i}R_{p_{i}0}+R_{p_{i}x^{i}}\right)+R_{x^{i}p_{i}}
=-\frac{n}{2}+\frac{n}{2}+\frac{n}{2}=\frac{n}{2}.\]

In summary, we have the following results:

\begin{proposition}
Metric $G$ has the Ricci Tensor
\begin{equation}
Ric(G)=(R_{ij})=\begin{pmatrix} -\frac{n}{2} & 0 & -\frac{n}{2}p_{j}\\
0 & 0 & \frac{1}{2}\delta_{ij}\\
-\frac{n}{2}p_{i} & \frac{1}{2}\delta_{ij} &
-\frac{n}{2}p_{i}p_{j}
\end{pmatrix},
\end{equation}
and the constant scalar curvature
\[
\mathcal{R}(G)=\frac{n}{2}.
\]

\end{proposition}
\par
To calculate the full curvature tensor recall that the curvature
transformation $R(X,Y):T(P)\rightarrow T(P)$ for vector fields
$X,Y$ is given by
\[R(X,Y)=\nabla_X\nabla_Y-\nabla_Y\nabla_X-\nabla_{[X,Y]}.\]
Note that R is antisymmetric in $X$ and $Y$.\par

Using formulas (6.4) for the covariant derivatives of basic
vectors of the frame (5.5) with respect to other basic vectors of
the same frame we get the action of the curvature transformation
$R(X,Y)$ with the couples $(X,Y)$ of basic vectors on all the
vectors of the frame (5.5).
\[R(\xi,P_i):  \xi\rightarrow\frac{1}{4}P_i;\hspace{.15 in} P_j\rightarrow
0;\hspace{.15 in}X_j\rightarrow -\frac{1}{4}\delta_{ij}\xi;\]
\[R(\xi,X_i):\xi\rightarrow\frac{1}{4}X_i;\hspace{.15 in}
P_j\rightarrow-\frac{1}{4}\delta_{ij}\xi ;\hspace{.15 in}
X_j\rightarrow 0\]
\[R(P_i,P_j):  \xi\rightarrow 0; \hspace{.15 in}P_k\rightarrow 0;\hspace{.15
in} X_k\rightarrow \frac{1}{4}(\delta_{ik}P_j-\delta_{jk}P_i)\]
\[R(P_i,X_j):  \xi\rightarrow 0; \hspace{.15 in}P_k\rightarrow
\frac{1}{4}\delta_{jk}P_i+\frac{1}{2}\delta_{ij}P_k;\hspace{.15
in} X_k\rightarrow
-\frac{1}{4}\delta_{ik}X_j-\frac{1}{2}\delta_{ij}X_k\]
\[R(X_i,X_j):  \xi\rightarrow 0; \hspace{.15 in}P_k\rightarrow
\frac{1}{4}(\delta_{ik}X_j-\delta_{jk}X_i); \hspace{.15
in}X_k\rightarrow 0.\] \par

By the antisymmetry properties of the curvature tensor R, this
determines the full tensor $R^{i}_{jkl}$.

\section{Sectional curvatures.}

Using the above calculations of the curvature transformations
$R(X,Y)$ we calculate, in the frame $(\xi ,P_{i}=\partial_{p_{i}}
,\partial_{x^{i}})$, the sectional curvatures of the nondegenerate
planes generated by couples of these vectors.\par

Given tangent vectors $A$ and $B$ spanning a surface in
$T_{m}(P)$, the sectional curvature they determine is given by the
formula (\cite{KN,Po})

\begin{equation}
R(m,A\wedge B)=\frac{G_{m}(R(A\wedge B)B,A))}{|A\wedge B|^2}
\end{equation}
where $|A\wedge B|^2=G(A,A)G(B,B)-G(A,B)G(B,A)$.

We calculate the sectional curvatures corresponding to the planes
generated by the following pairs of tangent vectors at a point
$m\in P$: $ (\xi, P_i), (\xi,
\partial_{x^i}), (P_i,P_j), (\partial_{x^i}, \partial_{x^j}),
(P_i,
\partial_{x^j}).$
\skip0.5cm

1. $(\xi,P_i)$
\[R(\xi,P_i)P_i=0\Rightarrow R(m,\xi \wedge P_{i})=0.\]
\skip0.5cm

2. $(\xi,\partial_{x^i})$
\[R(\xi,\partial_{x^i})\partial_{x^i}=R(\xi,X_i+p_i\xi)(X_i+p_i\xi)=R(\xi,X_i)(+p_i\xi)=+\frac{p_{i}}{4}X_i
   \]
Therefore,
\[ G(R(\xi,\partial_{x^i})\partial_{x^i}, \xi)=0\ \Rightarrow R(m,\xi \wedge
\partial_{x^i})=0.\]
\skip0.5cm

3. $(P_i,P_j)$
\[R(P_i,P_j)P_j=0\ \Rightarrow R(m,P_{i} \wedge P_{j})=0.\]
\skip0.5cm

4. $(\partial_{x^i}, \partial_{x^j})$
\[R(\partial_{x^i},
\partial_{x^j})\partial_{x^j}=R(X_i+p_i\xi,X_j+p_j\xi)(X_j+p_j\xi)=\]
\[=p_jR(X_i,\xi)((X_j+p_j\xi)+p_iR(\xi,X_j)(X_j+p_j\xi)=\]
\[=-\frac{p_j^2}{4}X_i+\frac{p_ip_j}{4}X_j.\]
Therefore, \[
G(R(\partial_{x^i},
\partial_{x^j})\partial_{x^j}, \partial_{x^i})=0\ \Rightarrow
R(m,\partial_{x^j} \wedge \partial_{x^j})=0.
\]\skip0.5cm

5. $(P_i,\partial_{x^j})$
\[R(P_i,\partial_{x^j})\partial_{x^j}=R(P_i,X_j+p_j\xi)(X_j+p_j\xi)=\]
\[=R(P_i,X_j)X_j+p_jR(P_i,\xi)X_j+p_j^2R(P_i,\xi)\xi=\]
\[=-\frac{3}{4}\delta_{ij}X_j+\frac{p_j}{4}\delta_{ij}\xi-\frac{p_j^2}{4}P_i
.\]
Therefore,
$G(R(P_i,\partial_{x^j})\partial_{x^j},P_i)=-\frac{3}{4}\delta_{ij}$.
\par
For the denominator we have $G(P_i,P_i)G(\partial_{x^j},
\partial_{x^j})-G(P_i,\partial_{x^j})G(\partial_{x^j},P_i)=-\delta_{ij}$.

Therefore, if $i=j$, the sectional curvature equals to $3/4$:
\[
R(m,P_i \wedge \partial_{x^i})=\frac{3}{4}.
\]
\par
When $i \ne j$, the sectional curvature is not defined since the
metric on the corresponding surface is degenerate (zero).
\begin{remark}
Notice that the only nonzero (positive !) sectional curvature in
this basis happens to be exactly at the planes corresponding to
the couples of conjugate variables $(p_{i},x^{i})$.
\end{remark}

\section{Killing vector fields, isotropy Lie algebra.}
In this section we calculate Killing vector fields for the metric
$G$ and, correspondingly, determine the Lie algebra of the
isotropy group of metric $G$.\par

Recall that such fields are defined by the condition
$\mathcal{L}_XG=0$. On the other hand for all $Y,Z$,
\[(\mathcal{L}_XG)(Y,Z)=XG(Y,Z)-G([X,Y],Z)-G(Y,[X,Z]).\]

From Section 7, formula (12), we have:
\[2G(X,\nabla_{Y}Z)=YG(Z,X)+ZG(X,Y)-XG(Y,Z)\]
\[+G(Y,[X,Z])+G([X,Y],Z)+G(X,[Y,Z])\]

\[=YG(Z,X)+ZG(Y,X)-(\mathcal{L}_XG)(Y,Z)+G(X,[Y,Z])\]\medskip

Since the connection is symmetric, $\mathcal{L}_XG=0$ if and only
if:
\begin{equation}
YG(X,Z)+ZG(X,Y)=G(X,\nabla_{Y}Z+\nabla_{Z}Y)
\end{equation}
for all $Y,Z$.  Notice that the condition $\mathcal{L}_XG(Y,Z)=0$
is linear by $Y,Z$.  Therefore, it is sufficient to ensure the
fulfillment of this condition for couples of vector fields $(Y,Z)$
from some frame on the manifold $P$.\par

\medskip

To determine Killing fields $X$, let $X=f\xi+g^kP_k+h^kX_k$.

Make the following substitutions into the formula (9.1): \[Y=\xi,
Z=\xi;\hspace{.15 in} Y=\xi,Z=P_i;\hspace{.15 in}
Y=\xi,Z=X_i\]\[Y=P_i,Z=P_j;\hspace{.15 in}Y=X_i,Z=X_j;\hspace{.15
in}Y=P_i,Z=X_j\]\medskip

Using the formulas (6.4) for covariant derivatives of basic
vectors along other basic vectors of basis $(\xi ,Z_i, P_j)$ we
get the system of equations for coefficients of the vector field
$X$:\bigskip

1.  $\xi f=0;$\medskip

2.  $h^i=P_i f+\xi h^i;$\medskip

3.  $g^i=-X_if-\xi g^i;$\medskip

4.  $P_ih^j+P_jh^i=0;$\medskip

5.  $X_ig^j+X_jg^i=0;$\medskip

6.  $P_ig^j+X_jh^i=0.$\bigskip

We begin by examining the case when $f=\theta(X)=0$.\bigskip

Since $\xi$ commutes with both $P_i$ and $X_j$, we may apply it to
equation 6 to obtain $P_i\xi g^j+X_j\xi h^i=0$.  By equations 2
and 3 this reduces to $P_ig^j-X_jh^i=0.$  Together with equation
6, this shows:\bigskip

7.  $P_ig^j=0, \hspace{.15 in} X_jh^i=0.$\bigskip

Now apply $X_k$ to equation 4.  Using the relation
$[X_k,P_l]=-\delta_{kl}\xi$ and equations 7, we see that
$0=-\delta_{ik}\xi h^j-\delta_{jk}\xi h^i$. By equation 2, this
gives $0=\delta_{ik}h^j+\delta_{jk}h^i$. Setting $i=j=k$, we
obtain $h^i=0$.  \bigskip

Similarly, if we apply $P_k$ to equation 5 and make the same
reductions, the result is $g^i=0$. Therefore, when $\theta(X)=0$,
we may conclude that $X=0$.\bigskip

Now consider the general case.  It is easy to see that $\xi$
satisfies the equations for a Killing vector field. Therefore, the
Lie bracket $[\xi, X]$ must also be a Killing vector field.  But by
equation 1, $\theta([\xi, X])=\xi f=0$.  So by the previous
calculations, $[\xi, X]=0$.  But this proves that $\xi h^i=0$, and
$\xi g^i=0$.\par   Hence equations 2 and 3 give explicit formulas
for $h^i$ and $g^i$ in terms of $f$ (a generating function).  If we
substitute these expressions into equations 4, 5, and 6, then 6 is
satisfied automatically, while 4 and 5 yield:
\[P_iP_jf=0, \hspace{.15 in} X_iX_jf=0.\]
Thus we see that $f$ must have the form:

\[f=a+b_ix^i+c^kp_k+d_i^kx^ip_k\]
with constant coefficients $a,b_i, c^k, d_k^i $.\par

From equations 2 and 3, we obtain $g^i$
and $h^i$:
\[g^i=b_i+d_i^kp_k;\]
\[h^i=-(c^i+d_k^ix^k).\]
As a result, the Lie algebra is generated by the vectors:
\begin{multline}
\xi;\ A_i= x^i\xi+P_i=\partial_{p_{i}}+x^{i}\partial_{x^{0}};\
B_j=p_j\xi-X_j=-\partial_{x^j},\\ Q_l^k=
x^kp_l\xi+p_lP_k-x^kX_l=p_l\partial_{p_k}-x^k\partial_{x^l}.
\end{multline}

\bigskip

These vector fields satisfy to the following nonzero commutator
relations:
\begin{equation}
[A^{i},B_{j}]=\delta_{ij}\xi,\hspace{.15 in}
[Q^{k}_{l},A_i]=-\delta_{il}A_k,\hspace{.15 in} [Q^{k}_{l},
B_{j}]=\delta_{kj}B_l,
\end{equation}
\begin{equation}
[Q^k_l,Q^r_s]=\delta_{ks}Q^r_l-\delta_{rl}Q^k_s
\end{equation}\medskip

As a result we get the following description of the Lie algebra
${\mathfrak iso}_{G}$ of the isometry group $Iso(G)$ of the metric
$G$

\begin{proposition}
The Lie algebra ${\mathfrak iso}_{G}$ of the isometry group
$Iso(G)$ of the metric $G$ is the Lie algebra ${\mathfrak
gl}(n,\mathbb{R})\times h_{n}$ being the semidirect product of
linear Lie algebra ${\mathfrak gl}(n,R)$ embedded into the
symplectic Lie algebra ${\mathfrak sp}(n,R)$ (with generators
$\{Q^k_l\}$) and of the Heisenberg Lie algebra ${\mathfrak h}_{n}$
with generators $\{\xi,A_i,B_j\}$ and commutative relations
(9.3-4).  All these vector fields are $\theta$-contact with the
contact Hamiltonians
$H_{\xi}=1,H_{A_{j}}=p_{j},H_{B_{i}}=x^{i},H_{Q^{k}_{l}}=
x^{k}p_{l}$ respectively (see \cite{AG}).
\end{proposition}
\vskip1cm

\section{Second fundamental form of Legendre surfaces.}
\vskip 0.5cm
\par
Here we will calculate the second fundamental form $II(X,Y)$ of a
Legendre submanifold $\Sigma_{\phi}$ at the points where metric
$g_{\phi}=G\vert_{\Sigma_{\phi}}$ is nondegenerate.
\medskip

Let $I$ be any subset of the indices from $1$ to $n$, and let $J$
be the complementary subset.  Consider a function $\phi(p_I, x^J)$
as in Sec.3 determining a Legendre submanifold. Throughout the
calculation, indices $i,i'$ will be assumed to be in $I$ and
indices $j,j'$ - in $J$. Any other indices will be assumed to run
through $1, \dots, n$. Also we will be using the following
conventions in notation:
\medskip

\[\phi_i=\partial_{p_i}\phi, \phi_j=\partial_{x^j}\phi\]

The Legendre submanifold $\Sigma_{\phi}$ is given, in these
notations by the equations (3.1):
\medskip

\[
\begin{cases}
x^0&=\phi - p_i\phi_{i},\\
p_j&=-\phi_{j},\\
 x^i&=\phi_i.
\end{cases}
\]
\medskip

The metric $G=2dp_k \odot dx^k+\theta^2$ restricted to this
submanifold $\Sigma_{\phi}$ is given by $g=2\phi_{ii'}dp_i \odot
dp_{i'}-2\phi_{jj'}dx^j \odot dx^{j'}$. Under the assumption that
$\Sigma_{\phi}$ is nondegenerate at a given point (and, therefore,
at some neighborhood of this point), the square matrices
$\phi_{ii'}, \phi_{jj'}$ must be nonsingular. Let their inverses
be given by $\phi^{ii'}, \phi^{jj'}$, respectively. Define, for
this section, the entries $\phi^{ij}$ to be $0$.
\medskip

Correspondingly, the tangent space $T(\Sigma_{\phi})$ to the
surface $\Sigma_{\phi}$ at each point is generated by the vectors:
\begin{multline}
Y_{i}=\phi_{*}(\partial_{p_{i}})=-p_{i'}\phi_{i'i}\partial_{x^0}+\partial_{p_i}-\phi_{ji}\partial_{p_j}+
\phi_{i'i}\partial_{x^{i'}}=
\partial_{p_i}-\phi_{ji}\partial_{p_j}+\phi_{i'i}(\partial_{x^{i'}}-p_{i'}\partial_{x^0})=\\
=P_{i}-\phi_{ij}P_{j}+\phi_{ii'}X_{i'}, ,\end{multline}
\begin{multline}
Y_{j}=\phi_{*}(\partial_{x^{j}})=(\phi_j-p_i\phi_{ij})\partial_{x^0}+\partial_{x^j}-\phi_{j'j}\partial_{p_{j'}}+\phi_{ij}\partial_{x^i}=
(\partial_{x^j}+\phi_{j}\partial_{x^{0}})-\phi_{j'j}\partial_{p_{j'}}+\phi_{ij}(\partial_{x^i}-p_i
\partial_{x^0})=\\ = X_{j}+\phi_{ij}X_{i}-\phi_{jj'}P_{j'}.\end{multline}
\medskip

We may simplify these expressions by defining vectors $V_k$ and
$W_k$ by the following rules.
\medskip

$V_k$ is equal to $P_k=\partial_{p_k}$ if $k$ belongs to $I$, and
is equal to $X_k=-p_k\partial_{x^0}+\partial_{x^k}$ if $k$ belongs
to $J$.
\medskip

$W_k$ equals $X_k$ when $k\epsilon I$ and equals $-P_k$ for
$k\epsilon J$.
\medskip

Then we see that the vectors above are given by one expression:
\begin{equation} Y_k=V_k+\phi_{kl}W_l.\end{equation}
\medskip
We also have the following useful relations for the scalar
products of these functions
\[
G(V_{k},V_{l})=G(W_{k},W_{l})=0,\ G(V_{k},W_{l})=\begin{cases}
\delta_{ks}, &\text{if $k,s\in I$;}\\
-\delta_{ks},&\text{if $k,s\in J$;}\\
0 ,&\text{otherwise.}
\end{cases}
\]

It is then easy to check that the following vectors are orthogonal
and complemental to ${Y_k}$:
\begin{equation} \xi,\
Z_k=W_k-\frac{1}{2}\phi^{kl}Y_l.\end{equation}

\medskip
To prove this we notice that $\xi $ is obviously orthogonal
$Y_{k}$ being orthogonal to $P_{i},X_{j}.$ Presenting  an
arbitrary vector in $D$ orthogonal to $Y_{k}$ in the form
$Z=a_{l}V_{l}+b_{k}W_{k}$ we find, using conditions of
orthogonality, relations between coefficients
\[
b_{i}=-\phi_{ii'}a_{i'}+\phi_{ij}a_{j},\
b_{j}=\phi_{ji}a_{i}-\phi_{jj'}a_{j'},
\]
and from this - the basic orthogonal vectors
\[Z_{i}=Y_{i}-2\phi_{ii'}W_{i'},\ Z_{j}=Y_{j}-2\phi_{jj'}W_{j'}.\]
Applying $\frac{1}{2}\phi^{il}$ to the right side of $Z_{i}$ (and
using the agreement that $\phi^{ij}=0$) we get from $Z_{i}$ the
other basic vectors $\frac{1}{2}\phi^{ii'}Y_{i'}-W_{i}$. Changing
sign and renaming these vectors $Z_{i}$ we get the second set of
orthogonal vectors in (10.4).  In the same way we get $Z_{j}$ of
the form (10.4).\par

In section 6, we computed the connection in the basis ${\xi, P_k,
X_l}$ (see (6.4)).  Referring to these calculations, it is simple
to see the following relations:
\[\nabla_{V_k}V_l=0, \nabla_{W_k}W_l=0;\]
\[\nabla_{V_k}W_l=-\nabla_{W_k}V_l=-\frac{1}{2}\delta_{kl}\xi .\]

Now we can easily calculate covariant derivatives of tangent
vector fields $Y_{k}$ with respect to $Y_{s}$:
\[\nabla_{Y_k}Y_l=\nabla_{(V_k+\phi_{kr}W_r})(V_l+\phi_{ls}W_s)=\]
\[=V_k(\phi_{ls})W_s+\phi_{ls}\nabla_{V_k}W_s+\phi_{kr}\nabla_{W_r}(V_l+\phi_{ls}W_s)=\]
\[=V_k(\phi_{ls})W_s-\frac{1}{2}\phi_{ls}\delta_{ks}\xi+\frac{1}{2}\phi_{kr}\delta_{rl}\xi+
\phi_{kr}W_r(\phi_{ls})W_s=\]
\[=V_k(\phi_{ls})W_s+\phi_{kr}W_r(\phi_{ls})W_s=\]
\[=Y_k(\phi_{ls})W_s=\phi_{lks}W_s.\]
\medskip

Expressed in the $Y_k, Z_l$ basis, the result is:
\[\nabla_{Y_k}Y_l=\phi_{lks}(Z_s+\frac{1}{2}\phi^{sr}Y_r)=\]
\[=\frac{1}{2}\phi^{sr}\phi_{lks}Y_r+\phi_{lks}Z_s.\]
\medskip

The first term in the last line represent the covariant derivative
$\nabla^\Sigma_{Y_k}Y_l$ on the submanifold $\Sigma_{\phi}$ with
respect to the induced (thermodynamical) metric while the second
term represents the second fundamental form of $\Sigma_{\phi} $
with respect to the Mrugala metric.\par

Namely we have proved the following
\begin{proposition}
Let a Legendre submanifold $\Sigma_{\phi}$ of the contact manifold
$(P,\theta)$ be defined by the equations
\[
\begin{cases}
x^0=\phi - \sum_{i\in I}p_i\phi_{i},\\
p_j=-\phi_{j},\\
x^i=\phi_i, \end{cases}, \ i\in I,j\in J. \]

Then the second fundamental form $II(X,Y)$ of the submanifold
$\Sigma_{\phi}$ is given by the expression
\[
II(Y_{k},Y_{l})=\phi_{lks}Z_s,
\]
where vector fields
\[Y_k=V_k+\phi_{kl}W_l\]
form the basis of tangent bundle of the submanifold $\Sigma_{\phi}$
(see above), vector fields $Z_k=W_k-\frac{1}{2}\phi^{kl}Y_l$ form
the basis of the orthogonal bundle of the submanifold
$\Sigma_{\phi}$ in the distribution $D$ and
$\phi_{lks}=Y_k(\phi_{ls})$.
\end{proposition}
\begin{example}
Consider the special case where $I= \emptyset$, so $\phi
=\phi(x^{j}),\ j=1,\ldots ,n$ and the tangent bundle to the surface
$\Sigma_{\phi}$ is generated by the tangent vectors
$Y_{j}=X_{j}-\phi_{jj'}P_{j'}$.  The Normal (orthogonal) subspace of
the tangent space $T_{m}(P)$ at the points of $\Sigma_{\phi}$ is
generated by
\[
Z_{j}=W_{j}-\frac{1}{2}\phi^{jj'}Y_{j'}=-P_{j}-\frac{1}{2}\phi^{jj'}(X_{j'}-\phi_{j'j''}P_{j''})=
\frac{1}{2}(P_{j}-\phi^{jj'}X_{j'}).
\]
We now calculate the coefficients $\phi_{lks}$ of the second
fundamental form ($k=j$ in this case)
\[
\phi_{ljs}=Y_{j}\phi_{ls}=(X_{j}-\phi_{jj'}P_{j'})\phi_{,x^{l}x^{s}}=\partial_{x^{j}}\phi_{,x^{l}x^{s}}.
\]
As a result, the second fundamental form of the surface
$\Sigma_{\phi}$ has the form
\begin{equation}
II(Y_{j},Y_{j'})=\phi_{,x^{j}x^{j'}x^{j''}}Z_{j''}=\eta_{\phi
,x^{j}}Z_{j},
\end{equation}
carrying information about all the third derivatives of a
thermodynamical potential $\phi $, or, equivalently, of the first
derivatives of the thermodynamical metric $\eta_{\phi}$.\par
Second fundamental form (10.5) of a surface $\Sigma_{\phi}$ is
zero iff the metric $\eta_{\phi}$ {\bf is constant with respect to
the variables $x^i$}. Only in such case (quite improbable in real
TD systems) metric $\eta_{\phi}$ is flat and submanifold
$\Sigma_{\phi}$ is {\bf totally geodesic} in $P$.
\end{example}
\vskip1cm

\section{Constitutive hypersurface.}
\vskip0.5cm

By the reasons of dimensions, the fundamental
thermodynamical constitutive equation (law) $x^{0}=\phi (x^{i})$ of
any material is {\bf homogeneous of order one}, i.e. the following
condition is fulfilled
\begin{equation}
\phi(\lambda x^{i})=\lambda \phi(x^{i})
\end{equation}
for all $\lambda \ne 0$ (see\cite{CA}). In other words, the action
of the one-parameter group $\Lambda $ of transformations
\begin{equation}
D_{\lambda}:(x^{0},p_{l},x^{i})\rightarrow (\lambda
x^{0},p_{l},\lambda x^{i}),\ \lambda \in R^{*}
\end{equation}
leaves the constitutive Legendre surface $\Sigma_{\phi}$ of a real
material invariant.
\par
As a result the surface $\Sigma_{\phi}$ lays in the canonical
quadric (hyperbolic paraboloid)
\begin{equation}
\mathcal{C}=\{ m=(x^{0},p_{l},x^{i})\vert
x^{0}+\sum_{i=1}^{n}p_{l}x^{l}=0\}.
\end{equation}
The intersection of the contact distribution $D$ with the fibers of
the tangent bundle $T(\mathcal{C})$ determines in $T(\mathcal{C})$
the subbundle $D_{\mathcal{C}}$. \par

Along the hypersurface $\mathcal{C}$ one has
\begin{equation}
0=dx^{0}+\sum_{i=1}^{n}p_{i}dx^{i}+\sum_{i}x^{i}dp_{i}=\theta
+\sum_{i}x^{i}dp_{i}.
\end{equation}
as a result, on the distribution $D_{\mathcal{C}}$ we have
\begin{equation}
\sum_{i}x^{i}dp_{i}=0
\end{equation}
which represents the abstract Gibbs-Duhem equation, \cite{CA}.\par
Using these relations it is easy to see that the subbundle
$D_{\mathcal{C}}$ contains and is generated by the
following vector fields
\begin{equation}
D_{\mathcal{C}}=<X_{l}=\partial_{x^{l}}+x^{l}\partial_{x^{0}},\
P_{ij}=x^{j}\partial_{p_{i}}-x^{i}\partial_{p_{j}}>
\end{equation}
at all points of $\mathcal C$ except the point of the plane ${\mathcal X}=\{
x^{i}=0,\ i=0,1,\ldots ,n\}$ where $D=T({\mathcal C})$.  To see
this we recall that $D$ is generated by the vector fields
$X_{l}=\partial_{x^{l}}+x^{l}\partial_{x^{0}},
P_{l}=\partial_{p_{l}}$.  The tangent space to the quadric $\mathcal C$
is formed by vectors satisfying the condition
\[
\theta +\sum_{i}x^{i}dp_{i}=0,
\]
Thus, the intersection $D\cap T({\mathcal C})$ is formed by vectors
from $D$ satisfying the Gibbs-Duhem equation. This is true for vectors
$X_{i}$ at all points.  Now, consider the partition of the index interval
$\{1,\ldots
,n\}=I\cup J$ and denote by ${\mathcal C}_{I}$ the set of points
\[
{\mathcal C}_{I}=\{ (x,p)\in {\mathcal C}\vert x^{i}=0,i\in
I;x^{j}\ne 0,j\in J\}.
\]
Let now $J\ne \varnothing $.  Then, for $m\in {\mathcal C}_{I}$
vectors $\partial_{p_{l}}, l\in I$ belong to $T_{m}({\mathcal C})$
as do all the vectors
$P_{ij}=x^{j}\partial_{p_{i}}-x^{i}\partial_{p_{j}},\ i,j\in J.$
Choose index $k\in J$ and consider vectors
$x^{l}\partial_{p_{k}}-x^{k}\partial_{p_{l}},\ l\in \{ J\backslash k
\}.$  These vectors belong to $T_{m}({\mathcal C})$, and, together
with $\partial_{p_{l}}, l\in I$ and $
X_{l}=\partial_{x^{l}}+x^{l}\partial_{x^{0}}$, they form the $2n-1$
dimensional subspace of $T_{m}({\mathcal C})$ which is the
intersection $D\cap T({\mathcal C})$.

The only points where these arguments fail to work are the points of
${\mathcal X}={\mathcal C}_{[1,n]}$, where $J=\varnothing $.  In
this case, the Gibbs-Duhem condition is empty and
$D_{m}=T_{m}({\mathcal C}).$\par

Notice that physically it would mean that all the extensive
variables of the system are zero - quite an improbable case.

\vskip1cm

\section{The Heisenberg Group as the thermodynamical phase space.}

In this section we establish isomorphism of the TPS $(P,\theta
,G)$ (with its contact structure and metric $G$) with the
Heisenberg Lie group $H_{n}$ with the right invariant contact
structure and right invariant indefinite metric.  This isomorphism
is locally suggested by the commutativity relations of vector
fields of the frame (5.5).\par

Recall that the Heisenberg group $H_{n}$ is the nilpotent Lie
group of $n\times n$ real matrices
\begin{equation}
g=g({\bar a}, {\bar b},c)=\begin{pmatrix} 1 & {\bar a} & c\\ 0 & I_{n} &
{\bar b}\\
0 & 0 & 1
\end{pmatrix}
\end{equation}
with the product
\begin{equation}
gg_{1}=\begin{pmatrix} 1 & {\bar a} & c\\ 0 & I_{n} & {\bar b}\\
0 & 0 & 1
\end{pmatrix} \cdot \begin{pmatrix} 1 & {\bar a}_{1} & c_{1}\\ 0 & I_{n} &
{\bar b}_{1}\\
0 & 0 & 1 \end{pmatrix} =
\begin{pmatrix} 1 & {\bar a}+{\bar
a}_{1} & c+c_{1}+
<{\bar a},{\bar b}_{1}>\\ 0 & I_{n} & {\bar b}+{\bar b}_{1}\\
0 & 0 & 1
\end{pmatrix},
\end{equation}
where $<{\bar a},{\bar b}_{1}>$ is the Euclidian scalar product of
two vectors from $R^n$ (see, for instance, \cite{OV}).\par

The Lie group $H_{n}$ is the central extension of the abelian
group $\mathbb{R}^{2n}_{{\bar a},{\bar b}}=\mathbb{R}^{n}_{\bar
a}\oplus \mathbb{R}^{n}_{{\bar b}}$ with the local parameters
${\bar a},{\bar b}$ respectively by the 1-dim abelian group
$\mathbb{R}_{c}$ with the local parameter $z$:
\[
1\rightarrow \mathbb{R}_{c}\rightarrow H_{n}\rightarrow
\mathbb{R}^{2n}_{{\bar a},{\bar b}}\rightarrow 1.
\]
\par
The Lie algebra ${\frak h}_{n}$ of the Heisenberg group is formed by
the matrices
\begin{equation}
X({\bar a},{\bar b},z)=\begin{pmatrix} 0 & {\bar a} & z\\ 0 & 0 & {\bar b}\\
0 & 0 & 0
\end{pmatrix}
\end{equation}
with the conventional matrix bracket as the Lie algebra operation.
The Lie algebra ${\mathfrak h}_{n}$ is mapped diffeomorphically onto
$H_{n}$ by the exponential mapping
\begin{equation}
exp(\begin{pmatrix} 0 & {\bar a} & z\\ 0 & 0 & {\bar b}\\
0 & 0 & 0
\end{pmatrix})=\begin{pmatrix} 1 & {\bar a} & c=z+\frac{1}{2}<{\bar a},{\bar
b}>\\ 0 & I_{n} & {\bar b}\\
0 & 0 & 1
\end{pmatrix}.
\end{equation}
\par
Lie algebra ${\frak h}_{n}$ is the central extension
\[
0\rightarrow \mathbb{R}_{z}\rightarrow {\frak h}_{n} \rightarrow
\mathbb{R}^{2n}_{{\bar a},{\bar b}}\rightarrow 0
\]
defined by the 2-cocycle $\omega_{+}:\mathbb{R}^{2n}\bigwedge
\mathbb{R}^{2n}\rightarrow \mathbb{R}$
\[
\omega_{+}(X({\bar a},{\bar b},z),X({\bar a}',{\bar
b}',z'))=\frac{1}{2}(<{\bar a},{\bar b}'>-<{\bar a}',{\bar b}>),
\]
of the canonical symplectic form in  $\mathbb{R}^{2n}$ see
(\cite{OV}).\par

We construct the diffeomorphic mapping
\[ \chi: H_{n}\Longleftrightarrow P
\]
by requiring
\begin{equation}
\chi :\ g=\begin{pmatrix} 1 & {\bar x} & x^{0}\\ 0 & I_{n} & {\bar p}\\
0 & 0 & 1
\end{pmatrix} \rightarrow m=\begin{pmatrix} -x^{0}\\{\bar p}\\{\bar
x}\end{pmatrix}.
\end{equation}
The Action of the group $H_{n}$ on itself by left translation:
$L_{g}: g_{1}\rightarrow gg_{1}$ defines the corresponding left
action of $H_{n}$ on the space $P$ as $T_{g}:m\rightarrow \chi
(L_{g} \chi^{-1}(m)))$, or
\begin{equation}
T_{g}\begin{pmatrix} x^{0}\\{\bar p}\\{\bar x}\end{pmatrix} =\chi
\left( \begin{pmatrix} 1 & {\bar a} & z\\ 0 & I_{n} & {\bar b}\\
0 & 0 & 1
\end{pmatrix}\cdot \begin{pmatrix} 1 & {\bar x} & -x^{0}\\ 0 & I_{n} & {\bar
p}\\
0 & 0 & 1
\end{pmatrix}   \right) =\begin{pmatrix} x^{0}-c-{\bar a}\cdot {\bar
p}\\{\bar p}+{\bar b}\\{\bar x}+{\bar
a}\end{pmatrix}.
\end{equation}

\par
Let us find, in these terms, generators of the left action of the
basic one-parameter subgroups of the group $H_{n}$ corresponding
to the elements $A_{i}=X(e_{i},0,0),\
B_{j}=X(0,e_{i},0),Z=X(0,0,1)$ of the Lie algebra ${\frak h}_{n}$
on the vectors of $P$.  For any element $X\in {\frak h}_{n}$
denote by $\xi_{X}$ (respectively by $\eta_{X}$ the right
invariant (respectively left invariant) vector field on $H_{n}$
generated by the left (respectively right) translations by
$exp(tX)$.
\par
We have in coordinates $({\bar A},{\bar b},c)$
\[
\xi_{Z}=\partial_{c},\
\xi_{A_{i}}=\partial_{a_{i}}+b_{i}\partial_{c},\
\xi_{B_{j}}=\partial_{b_{j}}.
\]
\par

Applying diffeomorphism $\chi $ to these vector fields we get the
correspondence
\begin{equation}
\begin{pmatrix} 1 & 0 & t\\ 0 & I_{n} & 0\\
0 & 0 & 1
\end{pmatrix} \cdot \chi^{-1}(m)=\begin{pmatrix} 1 & {\bar x} & -x^{0}+t\\ 0
& I_{n} & {\bar p}\\
0 & 0 & 1
\end{pmatrix} \rightarrow \chi_{*m}(\xi_{Z})=-\partial_{x^{0}},
\end{equation}
\begin{equation}
\begin{pmatrix} 1 & 0,0,\ldots t_{i},0\ldots & 0\\ 0 & I_{n} & 0\\
0 & 0 & 1
\end{pmatrix} \cdot \chi^{-1}(m)=\begin{pmatrix} 1 & {\bar x}+t{\bar e}_{i}
& x^{0}+p_{i}t\\ 0 & I_{n} & {\bar p}\\
0 & 0 & 1
\end{pmatrix} \rightarrow
\chi_{*m}(\xi_{A_{i}})=\partial_{x^{i}}-p_{i}\partial_{x^{0}}=X_{i},
\end{equation}
\begin{equation}
\begin{pmatrix} 1 & 0 & 0\\ 0 & I_{n} & 0,0,\ldots t_{j},0\ldots\\
0 & 0 & 1
\end{pmatrix} \cdot \chi^{-1}(m)=\begin{pmatrix} 1 & {\bar x} & x^{0}\\ 0 &
I_{n} & {\bar p}+t{\bar f}_{j}\\
0 & 0 & 1
\end{pmatrix} \rightarrow \chi_{*m}(\xi_{B_{j}})=\partial_{p_{j}}=P_{j}.
\end{equation}

\par
The pullback of the contact form $\theta =dx^{0}+p_{l}dx^{l} $ from
$P$ to $H_{n}$ defines the 1-form $\theta_{H}$ on $H_{n}$
\[
\theta_{H}=\chi^{*}(\theta )=-dc+b_{i}da^{i}.
\]
\par
 Reeb vector field of this form is
\[
\xi_{H}=-\xi_{c}=-\partial_{c},
\]
and we have
\[
\chi_{*}(\xi_{H})=\xi
\]
for the Reeb vectors of contact manifolds $(H_{n},\theta_{H})$ and
$(P,\theta )$.\par

The kernel, $D_{H}$, of this 1-form (a distribution of codimension 1
on $H_{n}$) is, at each point $g$, generated by the values of vector
fields $\xi_{A_{i}},\xi_{B_{j}}$ of left translations, and is
therefore {\bf right invariant}.  As a result, distribution $D_{H}$
defines the {\bf right invariant contact structure on $H_{n}$}
(given as the kernel of the form $\theta_{H}$).  \par

Considering the {\bf right translations} on the group $H_{n}$,
corresponding to the 1-dim Lie subalgebras of ${\frak h}_{n}$ with
generators $Z,A_{i},B_{j}$ we find that their generators have the
form
\begin{equation}
\eta_{C}=\partial_{c},\ \eta_{A_{i}}=\partial_{a_{i}},\
\eta_{B_{j}}=\partial_{b_{j}}+a_{i}\partial_{c}
\end{equation}
and it is easy to check that the form $\theta $ {\bf is invariant
under the flow of these (therefore contact) vector fields}.  Thus,
the form $\theta_{H}$ is {\bf right invariant}. The diffeomorphism
$\chi $ send these vectors into
\[
\chi_{*}(\eta_{C})=\xi;\ \chi_{*}(\eta_{A_{i}})=\partial_{x^{i}}
;\ \chi_{*}(\eta_{B_{j}})=\partial_{p_{j}}+x^{j}\partial_{x^{0}}.
\]
Comparing this result with the description of the Killing vector
field of the metric $G$ we see that these vector fields form the
nilradical of the Lie algebra of the Killing vector fields of the
metric $G$.
\par

\begin{remark}
The distribution $D_{H}$ is the direct sum of two n-dimensional
distributions
\[
D=D_{A}\oplus D_{B}
\]
in obvious notations, distributions $D_{A},D_{B}$ are integrable
having as the basis at each point values of pairwise commuting
vector fields.  Denote by $A$ (respectively by $B$) the abelian
subgroup of $H_{n}$ of matrices of the form (12.1) with $c={\bar
b}=0$ (respectively $c={\bar a}=0$).  Then integral manifolds of
distribution $D_{A}$ (respectively of distribution $D_{B}$) are
orbits of the left translations by the subgroup $A$ (respectively,
by the Lie subgroup $B$).
\end{remark}

Metric $G$ is transferred under the diffeomorphism $\chi $ into
the metric $G_{H}$ on the Heisenberg group.  This metric is
constant in the right invariant (non-holonomic) frame
$(\xi_{C},\xi_{A_{i}},\xi_{B_{j}})$ and is, therefore, right
invariant by itself.\par

As a result we've proved the following
\begin{theorem}
The diffeomorphism $\chi $ defined by
\[
\chi :\ g=\begin{pmatrix} 1 & {\bar x} & x^{0}\\ 0 & I_{n} & {\bar p}\\
0 & 0 & 1
\end{pmatrix} \rightarrow m=\begin{pmatrix} -x^{0}\\{\bar p}\\{\bar
x}\end{pmatrix}.
\]
determines an isomorphism of
the "thermodynamical metric contact manifold" $(P,\theta ,G)$ with
the Heisenberg group $H_{n}$ endowed with the right invariant
contact from $\theta_{H}$ and the right invariant metric $G_{H}$
of signature $(n+1,n).$
\end{theorem}

\par
\begin{remark}
Recall (see \cite{GP}) that the automorphism  group $Aut(H_{n})$
of the Heisenberg group has, as its connected component of unity
the Lie group $Aut_{0}(H_{n})=Sp(n,\mathbb{R})\times A^{2n+1}$,
where $A^{2n+1}$ is the abelian group of dimension $2n+1$. This
group acts on the space of all right invariant contact 1-forms on
the group $H_{n}$. Since right invariant one-forms on $H_{n}$ are
defined by their values at the unit of the group $e\in H_{n}$ and
since the automorphism group $Aut(G)$ leaves $e\in G$ fixed it is
sufficient to study action of this group at the set of elements
$\nu \in {\frak h}_{n}^{*}.$  The following result for the left
invariant contact structures (with the sketch of the proof) was
sent to the authors in a letter by M. Goze. We reformulate this
result for the right invariant contact structures due to the
obvious duality between left and right translations.\par

Let $\{  X_{i},Y_{j},Z\}$ be a standard basis of ${\frak h}_{n}$
with the only nontrivial brackets being $[X_{i},Y_{i}]=Z,\
i=1,\ldots ,n$.  Let $\{  \alpha_{i},\beta_{j},\omega\}$ be the
dual basis in ${\frak h}_{n}^{*}$.  Then, extending this basis to
the coframe of right invariant vector fields we get relations
\[
d\alpha_{i}=d\beta_{j}=0,\ d\omega
=\sum_{i}\alpha_{i}\wedge\beta_{i}.
\]
Then $\omega\wedge (d\omega)^{n}\ne 0$ and, therefore, one-form
$\omega $ defines the right invariant contact structure on
$H_{n}$.\par

\begin{proposition}
Let $\omega_{1}=a\omega
+\sum_{i}a_{i}\alpha_{i}+\sum_{j}b_{j}\beta_{j}$ be a contact form
in ${\frak h}_{n}^{*}$. The group $Aut(H_{n})$ acts transitively
on the set of right invariant contact structures with the isotropy
group of $\omega $ being the intersection of the group
$Aut(H_{n})$ with the group of $\omega $-conformally contact
diffeomorphisms of $H_{n}$.
\end{proposition}
It follows from this that the contact structure $\theta_{H}$ of the
thermodynamical phase space is the typical representative of the
$Aut(H_{n})$-conjugacy class of right invariant contact structures
on the Heisenberg Group defined by a choice of the canonical basis
$\{ X_{i},Y_{j},Z\}$ of the Lie algebra ${\frak h}_{n}$ and,
therefore, unique, up to an automorphism of the group $H_{n}$.
\end{remark}
\begin{remark}
After the isomorphism of the TPS $(P,\theta, G)$ with
$(H_{n},\theta_{H}, G_{H})$ is established, many properties of
metric $G$ can be obtained from the corresponding results for
invariant metrics on Lie groups. Further use of this isomorphism
for the study of thermodynamical systems will be a subject of
future work.
\end{remark}

\centerline{ \textbf{Part II.}}
\section{Symplectization of manifold $(P,\theta,G)$.}

Let $\tP$ be the standard (2n+2)-dim real vector space
$\mathbb{R}^{2n+2}$ with the coordinates $(p_{i},x^{j}),\
i,j=0,\ldots n$, endowed with the 1-form
\begin{equation}
{\tilde \theta }= \sum_{i=0}^{n}p_{i}dx^{i},
\end{equation}
and the standard symplectic structure
\begin{equation}
\omega =d\theta = \sum_{i}dp_{i}\wedge dx^{i}.
\end{equation}
We consider the {\bf embedding} of the space
$(P,\theta=dx^{0}+p_{l}dx^{l})$ into $\tP$
\begin{equation}
J: (x^{0},x^{i},p_{j})\rightarrow (x^{0},x^{i};p_{0}=1,p_{l},\
l=1,\ldots ,n)
\end{equation}
as the affine subspace $p_{0}=1.$\par It is easy to see that
\begin{proposition}
\begin{enumerate}
\item  The pullback by $J$ of the 1-form $\tilde \theta$ coincides with
the contact form $\theta$.
\[
J^{*}({\tilde \theta})=\theta .
\]
\item The symplectic manifold $(\tP =\{(p,x)\in
\mathbb{R}^{2n+2}\vert p_{0}>0\},\omega )$ is the standard {\bf
symplectization} of $(P,\theta)$ (see \cite{A,LM}) and $J$ is the
section of the symplectization bundle $\pi: \tP \rightarrow P.$

\item The symmetrical tensor
\begin{equation}
{\tilde G}=(\tG_{ij})=
\begin{pmatrix}
  0_{n+1\times n+1}  &  I_{n+1}\\
  I_{n+1} & p_{i}p_{j}
\end{pmatrix}
\end{equation}
determines in $\tP$ the pseudo-Riemannian metric of signature
$(n+1,n+1)$.
\item The restriction of metric $\tG$ to the image of
the embedding $J$ coincides with the metric $G$
\[
J^{*}\tG=G.
\]
\item There is a bijection between the Legendre submanifolds of
the contact manifold $(P,\theta )$ and the {\it homogeneous} (under
the action $(x^{i},p_{j})\rightarrow (x^{0}, \lambda
p_{0},x^{1},p_{1},\ldots, x^{n},p_{n})$ of $R^+$ on the manifold
$\tP$) Lagrange submanifolds of the symplectic manifold $(\tP
,\omega))$. This correspondence is defined by the intersection of a
homogenous Lagrangian submanifold $\tilde K$ with the image of the
embedding $J$ and by the action of the dilatation group on the image
of a Legendre submanifold $K\subset P$ under the embedding $J$.
\end{enumerate}
\end{proposition}
\begin{proof}
Almost all the statements of this Proposition follows simply from
the construction or are known (\cite{AG,LM}). The determinant of
the matrix (13.4) of the metric $\tilde G$ is equal to
$(-1)^{n+1}$ which proves its nondegeneracy.
\end{proof}
\vskip0.5cm

\section{Canonical frame $(\tX^{i},\tP_{j})$.}

Here we introduce the non-holonomic frame $(\tX^{i},\tP_{j})$ in
the open subset $\tP^{*}=\{ (p,x)\in \tP \vert \prod_{k}p_{k}\ne
0\} $ of the symplectic manifold $\tP$ with respect to which the
metric $\tG$ has the standard (constant) form.\par

We take
\begin{equation}
\tP_{i}=p_{i}\partial_{p_{i}};\ L_{k}=p_{k}^{-1}\partial_{x^{k}};\
\tX_{j}=p_{k}^{-1}\partial_{x^{k}}- {\hat P}=L_{j}
- {\hat P},
\end{equation}
where
\[
{\hat
P}=\frac{1}{2}\sum_{s=0}^{n}\tP_{s}=\frac{1}{2}\sum_{i=0}^{n}p_{i}\partial_{p_{i}}
\]
is the generator of homogeneous dilatation in the $p$-directions.
\par

Then we have the following commutator relations between
the introduced vector fields.
\begin{multline}
[ \tP_{i}, \tP_{j}]=[L_{i},L_{j}]=0;\ [\tP_{i},\
L_{j}]=-\delta_{ij}L_{j};\ [ L_{j}, {\hat P} ]=\frac{1}{2}L_{j};\\
[\tP_{i}, \tX_{j}]=-\delta_{ij}(\tX_{j}+{\hat P});
[\tX_{i},\tX_{j}]=\frac{1}{2}(\tX_{j}-\tX_{i}).
\end{multline}
These relations shows, in particular,  that the couples of vector
fields $\tP_{i},L_{i}$ form the 2-dim solvable Lie algebras
$\alpha_{i}$ of vector fields commuting between themselves (see
below, Sec.19).\par

Scalar products of the introduced vector fields are calculated as
follows
\begin{multline}
\tG(\tP_{i},\tP_{j})=0;\ G(L_{i},L_{j})=1;\ G(P_{i},L_{j})=\delta_{ij};\\
\tG(\tP_{i}, \tX_{j})=\tG(p_{i}\partial_{p_{i}},
p^{-1}_{j}\partial_{x^{j}}-{\hat
P})=p_{i}p_{j}^{-1}\tG(\partial_{p_{i}},\partial_{x^{j}})=p_{i}p_{j}^{-1}\delta_{ij}=\delta_{ij};\\
\tG(\tX_{i},{\hat
P})=\frac{1}{2}\sum_{s}\tG(\tX_{i},\tP_{s})=\frac{1}{2}\sum_{s}\delta_{is}=\frac{1}{2};\\
\tG(\tX_{i},\tX_{j})= \tG(p^{-1}_{i}\partial_{x^{i}}-{\hat
P},p^{-1}_{j}\partial_{x^{j}}-{\hat
P})=p_{i}^{-1}p_{j}^{-1}\tG(\partial_{x^{i}},\partial_{x^{j}})-p_{i}^{-1}\tG(\partial_{x^{i}},{\hat
P})-\\ p_{j}^{-1}\tG(\partial_{x^{j}},{\hat P})+\tG({\hat P},{\hat
P})=
1-\frac{1}{2}p_{i}^{-1}\sum_{s}p_{s}\tG(\partial_{x^{i}},\partial_{p_{s}})
-\frac{1}{2}p_{j}^{-1}\sum_{s}p_{s}\tG(\partial_{x^{j}},\partial_{p_{s}})+0=\\
1-\frac{1}{2}p_{i}^{-1}\sum_{s}p_{s}\delta_{is}-\frac{1}{2}p_{j}^{-1}\sum_{s}p_{s}\delta_{js}=1-1=0.
\end{multline}
As a result in the basis $(\tP_{i}, \tX_{j})$, the matrix of
metric $\tG$ has in $\tP^{*}$ the following canonical form
\begin{equation}
(\tG_{ij})=\begin{pmatrix} 0_{n+1} & I_{n+1}\\ I_{n+1} & 0_{n+1 }
\end{pmatrix} .
\end{equation}
\vskip0.5cm
\par
The positive and negative distributions of metric $\tG$ are
\begin{equation}
T^{+}=<\xi_{i}=\frac{1}{\sqrt{2}}(\tP_{i}+\tX_{i})>,\
T^{-}=<\eta_{i}=\frac{1}{\sqrt{2}}(\tP_{i}-\tX_{i})>.
\end{equation}
\par
The zero cone is given in the frame $\tP_{i},\tX_{j}$ for
$X=f_{i}\tP_{i}+g_{j}\tX_{j}$
\begin{equation}
\tG(X,X)=0 \Leftrightarrow \sum_{i}f_{i}g_{i}=0.
\end{equation}

\section{Levi-Civita connection of metric $\tG$.}
To calculate curvature of the metric $\tG$ we start with the
combinations $\{ ij,k\}=\tG_{ik,j}+\tG_{jk,i}-\tG_{ij,k}$ It is
easy to see that the only nonzero combinations are
\begin{equation}
\{ x^{i}x^{j},p_{l}\}=-(\delta^{i}_{l}p_{j}+\delta^{j}_{l}p_{i});\
\{ x^{i}p_{j},x^{s}\}=\{
p_{j}x^{i},x^{s}\}=(\delta^{i}_{j}p_{s}+\delta^{s}_{j}p_{i}).
\end{equation}
Using these combinations we calculate the Christoffel coefficients
$\Gamma^{i}_{jk}=\frac{1}{2}G^{is}\{ jk,s\}.$ We notice that to be
nonzero, the Christoffel coefficient $\Gamma^{\mu }_{\nu \xi}$
should have at least two $x^{i}$ between the three indices in one of
the terms $\{ jk,s\}$. Using this it is easy to see that the only
nonzero Christoffel coefficients are
\begin{equation}
\Gamma^{x^{i}}_{x^{j}x^{k}}=-\frac{1}{2}(\delta^{i}_{j}p_{k}+\delta^{i}_{k}p_{j});\
  \Gamma^{p_{i}}_{x^{j}x^{k}}=p_{i}p_{j}p_{k};\
\Gamma^{p_{i}}_{x^{j}p_{k}}=\frac{1}{2}(\delta^{j}_{k}p_{i}+\delta^{i}_{k}p_{j}).
\end{equation}

Next, we calculate the Ricci Tensor.
$R_{ij}=\Gamma^{k}_{ji,k}-\Gamma^{k}_{js}\Gamma^{s}_{ki}$. We have
\[
R_{p_{i}p_{j}}=\Gamma^{k}_{p_{i}p_{j},k}-\Gamma^{k}_{p_{j}s}\Gamma^{s}_{kp_{i}}=
0-\Gamma^{p_{s}}_{p_{j}x^{l}}\Gamma^{x^{s}}_{p_{s}p_{i}}-
\Gamma^{x^{s}}_{p_{j}x^{l}}\Gamma^{x^{l}}_{x^{s}p_{i}}=0,
\]
Since
$\Gamma^{p_{s}}_{x^{i}x^{j},p_{s}}=\sum_{s}\partial_{p_{s}}(p_{i}p_{j}p_{s})=(n+1)p_{i}p_{j}+
\sum_{s}(\delta_{s}^{i}p_{j}p_{s}+\delta^{j}_{s}p_{i}p_{s})=(n+1)p_{i}p_{j}+p_{i}p_{j}+p_{i}p_{j}=(n+3)p_{i}p_{j},$
we have
\begin{multline}
R_{x^{i}x^{j}}=\Gamma^{k}_{x^{i}x^{j},k}-\Gamma^{k}_{x^{j}s}\Gamma^{s}_{kx^{i}}=
\Gamma^{p_{s}}_{x^{i}x^{j},p_{s}}-\Gamma^{x^{l}}_{x^{j}x^{s}}\Gamma^{x^{s}}_{x^{l}x^{i}}-
\Gamma^{x^{l}}_{x^{j}p_{s}}\Gamma^{p_{s}}_{x^{l}x^{i}}-
\Gamma^{p_{l}}_{x^{j}x^{s}}\Gamma^{x^{s}}_{p_{l}x^{i}}-
\Gamma^{p_{l}}_{x^{j}p_{s}}\Gamma^{p_{s}}_{p_{l}x^{i}}=\\
=(n+3)p_{i}p_{j}-\frac{1}{4}(\delta^{l}_{j}p_{s}+\delta^{l}_{s}p_{j})
(\delta^{s}_{l}p_{i}+\delta^{s}_{i}p_{l})-
\frac{1}{4}(\delta^{j}_{s}p_{l}+\delta^{l}_{s}p_{j})
(\delta^{i}_{l}p_{s}+\delta^{s}_{l}p_{i})=\\ = (n+3)p_{i}p_{j}-
\frac{1}{4}[p_{i}p_{j}+p_{i}p_{j}+(n+1)p_{i}p_{j}+p_{i}p_{j}]-
\frac{1}{4}[p_{i}p_{j}+p_{i}p_{j}+(n+1)p_{i}p_{j}+p_{i}p_{j}]=\\
=(n+3)p_{i}p_{j}-
\frac{1}{2}(n+4)p_{i}p_{j}=\frac{n+2}{2}p_{i}p_{j}.
\end{multline}
Finally, since
$\Gamma^{p_{s}}_{x^{j}p_{i},p_{s}}=\frac{1}{2}\sum_{s}\partial_{p_{s}}(\delta_{i}^{j}p_{s}+\delta^{s}_{i}p_{j})=
\frac{1}{2}((n+1)\delta_{i}^{j}+\delta^{i}_{j})=\frac{n+2}{2}\delta^{i}_{j},$
we have
\begin{multline}
R_{p_{i}x^{j}}=\Gamma^{k}_{x^{j}p_{i},k}-\Gamma^{k}_{x^{j}s}\Gamma^{s}_{kp_{i}}=
\Gamma^{p_{s}}_{x^{j}p_{i},p_{s}}-(
\Gamma^{p_{l}}_{x^{j}x^{k}}\Gamma^{x^{k}}_{p_{l}p_{i}}+
\Gamma^{p_{l}}_{x^{j}p_{s}}\Gamma^{p_{s}}_{p_{l}p_{i}}+
\Gamma^{x^{l}}_{x^{j}x^{s}}\Gamma^{x^{s}}_{x^{l}p_{i}}+
\Gamma^{x^{l}}_{x^{j}p_{s}}\Gamma^{p_{s}}_{x^{l}x^{i}})=\\
=\frac{n+2}{2}\delta^{j}_{i}+0-0-0-0-0=\frac{n+2}{2}\delta^{j}_{i}.
\end{multline}
Thus,
\begin{equation}
Ric(\tG)=\begin{pmatrix} 0_{(n+1)\times (n+1)} &
\frac{n+2}{2}I_{n+1}\\
\frac{n+2}{2}I_{n+1} & \frac{n+2}{2}p_{i}p_{j}
\end{pmatrix}= \frac{n+2}{2}\tG.
\end{equation}
Thus, metric $\tG$ is {\bf pseudo-Riemannian Einstein metric},see
\cite{ON}.\par Lifting an index of $R_{ij}$ with
\[
(G^{ij})=\begin{pmatrix} -p_{i}p_{j} &
I_{n+1}\\
I_{n+1} & 0_{(n+1)\times (n+1)}
\end{pmatrix}
\]
we get $\tG^{ik}R_{kj}=\frac{n+2}{2}\delta^{i}_{j}$.  Taking the
trace, we find the scalar curvature to be
\begin{equation}
R(\tG)=Tr(\tG^{-1}\frac{n+2}{2}\tG)=(n+1)(n+2).
\end{equation}
Thus, we have
\begin{proposition}
Metric $\tG$ is the indefinite Einstein metric of scalar curvature
$R(\tG)=(n+1)(n+2)$.
\end{proposition}

\section{Killing vector fields of metric $\tG$.}
It is natural to find the form of Killing vector fields of metric
$\tG$ - infinitesimal isometries of $\tG$. The details of the
calculations are given in the Appendix. Here we formulate the
final result.\par

\begin{theorem}
The Lie algebra ${\mathfrak iso}_{\tG}\simeq{\mathfrak
{sl}}(n+2,\mathbb{R})$ of Killing vector fields of the metric
$\tG$ is (as the vector space) the linear sum
\[
{\mathfrak iso}_{\tG}= {\mathfrak q}\oplus {\mathfrak d}\oplus
{\mathfrak x }
\]
of Lie subalgebras
\[
1)\ {\mathfrak
q}=<Q^{i}_{j}=x^{i}\partial_{x^{j}}-p_{j}\partial_{p_{i}}>,
\]
with the commutator relations
\[
[Q^{i}_{j},Q^{p}_{k}]=\delta^{p}_{j}Q^{i}_{k}-\delta^{i}_{k}Q^{p}_{j},
\]
Subalgebra $\mathfrak q$ is isomorphic thus, to ${\mathfrak
gl}(n+1,\mathbb{R})$,\par

the abelian subalgebra
\[
2)\ {\mathfrak d}=<D^{i}=
\frac{x^{i}}{2}Q+(1-\frac{1}{2}(x^{l}p_{l}))\partial_{p_{i}},>
\]
where
$Q=\sum_{i}Q^{i}_{i}=\sum_{i}(x^{i}\partial_{x^{i}}-p_{i}\partial_{p_{i}})$
is the generator of hyperbolic rotation $H_{t}:(p,x)\rightarrow
(e^t p, e^{-t} q)$.
\par
and abelian subalgebra
\[
3)\ {\mathfrak x }=<X_{s}=\frac{\partial}{\partial x^{s}}>.
\]
\par
Generators $Q^{i}_{j},X_{j},D^{j}$ satisfy to the following
commutator relations
\begin{equation}
[Q^{i}_{j},X_{s}]=-\delta^{i}_{s}X_{j};\
[Q^{i}_{j},D^{s}]=\delta^{s}_{j}D^{i};\
[X_{s},D^{i}]=\frac{1}{2}Q^{i}_{s} +\frac{1}{2}\delta^{i}_{s}Q.
\end{equation}
\par
Vector fields $Q^{i}_{j},X_{i},D^{j}$ are Hamiltonian with
Hamiltonian functions
\[
H_{Q^{i}_{j}}=-x^i p_j;\ H_{X_k}=-p_k;\ H_{D^s}=x^s
(1-\frac{<{\bar x},{\bar p}>}{2}).
\]
\end{theorem}

\section{Hypersurface $\mathcal{\tilde C}$.}

There exists a natural lift of the constitutive hyperquadric
$\mathcal{C}$ to the space $\tP$ as the homogeneous hypersurface
\begin{equation}
\mathcal{\tilde C}=\{ (p_{i},x^{i})\vert
x^{0}p_{0}+\sum_{i=1}^{i=n}p_{i}x^{i}=0
\}.
\end{equation}
The hyperquadric $\mathcal{C}$ is invariant under the action of
$\Lambda = \mathbb{R}^{*}$: $(x^{i},p_{i})\rightarrow (
x^{i},\lambda p_{i})$ as well as under the hyperbolic rotations
$(x^{i},p_{i})\rightarrow (\lambda x^{i},\lambda^{-1} p_{i}).$\par

Polarizing the coordinates $(x^{i},p_{j})$ - introducing new
coordinates $\xi_{i}=x^{i}+p_{i},\ \eta_{j}=x^{j}-p_{j}$ we rewrite
the equation of $\mathcal{\tilde C}$ as follows
\[
\sum_{i=0}^{n}(\xi_{i}^{2}-\eta_{i}^{2})=0.
\]
From this we see that the hyperquadric $\mathcal{\tilde C}$ is a
cone in the space $\tilde P$ of signature $(n+1,n+1)$.

\vskip1cm

\section{$(P,\phi ,\theta ,\xi )$ as the indefinite Sasakian manifold.}
Recall that an almost contact manifold $(M^{2n+1},\phi ,\xi,
\eta)$ is called {\bf Sasakian} if the {\bf almost complex
structure} on the manifold $M^{2n+1}\times R$ defined by
\begin{equation}
J(X,f\partial_{t})=(\phi (X)-f\xi ,\eta(X))\partial_{t})
\end{equation}
{\bf is integrable}.  Here $t$ is the coordinate on the factor $R$
of the product, $f\in C^{\infty}(M\times R)$, see \cite{BL}.\par

Sasakian manifolds are considered to be the natural
odd-dimensional analog of {\bf Kahler} manifolds,\cite{BL},
Chapter 6.\par

A necessary and sufficient condition for the integrability of an
almost complex structure is the vanishing of the {\bf Nijenhuis
tensor} of the (1,1)-tensor $J$:
\[
N_{J}(X,Y)=J^{2}[X,Y]+[JX,JY]-J[JX,Y]-J[X,JY].
\]
In the situation where $(M,\theta )$ is the {\bf contact manifold}
and the almost contact structure $(\phi ,\xi, \eta )$ is
associated with the contact structure (so that in particular, the
1-form $\eta =\theta $ is the same in both structures and $\xi $
is the Reeb vector of the contact structure) it is natural to
study the integrability of the almost complex structure defined by
(18.1) on the {\bf symplectization} of the manifold $(M,
\theta).$\par

In the case of the standard contact structure $(P,\theta )$, the
symplectization of $P$ is naturally embedded in the symplectic
vector space $(\tP,\omega )$ which can be considered as the
product manifold $P\times \mathbb{R}$.\par Even in the case of a
general contact manifold this seems to be a natural modification
of the definition of "normality" of an almost contact structure.

\begin{proposition}
The Nijenhuis tensor $N_{J}$ of the almost complex structure
defined by the formula (18.1) on the symplectization $\tP$ of the
manifold $(P,\phi ,\xi ,\theta)$ with $\phi$ given in (5.12) {\bf
is identically zero}. As a result, $(\tP , J)$ is a complex
manifold and $(P,\phi ,\xi ,\theta)$ is an {\bf "indefinite
Sasakian manifold"}.
\end{proposition}
\begin{proof}
With a slight abuse of notations we will use the coordinate $t$
instead of $p_{0}$ for the $n+2$-th coordinate in $\tP$.  By the
linearity of the condition $N_{J}(X,Y)=0$, it is sufficient to only
check vectors from some frame. We will use the frame $(\xi ,
X_{i},P_{j},\partial_{t})$ for the calculations.\par

Recall that
\[
\phi (X_{i}=P_{i},\phi (P_{i})=-X_{i},\ \phi(\xi
)=0.
\]
\par

With this we get
\[
J(X)=\phi(X)+\theta(X)\partial_{t},\ X\in T(P);\ J(\xi
)=\partial_{t},\ J(\partial_{t})=-\xi .
\]

Now we calculate
\[
N_{J}(\xi ,\partial_{t})=J^{2}[\xi ,\partial_{t}]+[J\xi
,J\partial_{t}]-J[J\xi ,\partial_{t}]-J[\xi
,J\partial_{t}]=0+[\partial_{t}, -\xi
]-J[\partial_{t},\partial_{t}]-J[\xi ,-\xi]=0,
\]

\[
N_{J}(X_{i},\partial_{t})=J^{2}[X_{i} ,\partial_{t}]+[\phi X_{i}
,-\xi ]-J[\phi X_{i} ,\partial_{t}]-J[X_i ,-\xi ]=-[P_{i},\xi
]-J[P_{i},\partial_{t}]-J0=0,
\]

\[
N_{J}(P_{i},\partial_{t})=J^{2}[P_{i} ,\partial_{t}]+[-X_{i} ,-\xi
]-J[- X_{i} ,\partial_{t}]-J[P_i ,-\xi ]=0,
\]

\[
N_{J}(X_{i},\xi )=J^{2}[X_{i},\xi
]+[P_{i},\partial_{t}]-J[P_{i},\xi ]-J[X_{i}, \partial_{t}]=0,
\]

\[
N_{J}(P_{i},\xi)=J^{2}[P_{i},\xi
]+[-X_{i},\partial_{t}]-J[-X_{i},\xi ]-J[P_{i}, \partial_{t}]=0,
\]
\[
N_{J}(X_{i},X_{j})=J^2
[X_{i},X_{j}]+[P_{i},P_{j}]-J[P_{i},X_{j}]-J[X_{i},P_{j}]=-J(-\delta_{ij}\xi
)-J(\delta_{ij}\xi )=0,
\]
\[
N_{J}(P_{i},P_{j})=J^2
[P_{i},P_{j}]+[-X_{i},-X_{j}]-J[-X_{i},P_{j}]-J[P_{i},-X_{j}]=-J(\delta_{ij}\xi
)+J(-\delta_{ij}\xi )=0,
\]
\[
N_{J}(X_{i},P_{j})=J^2
[X_{i},P_{j}]+[P_{i},-X_{j}]-J[P_{i},P_{j}]-J[X_{i},-X_{j}]=J^{2}[\delta_{ij}\xi
]+\delta_{ij}\xi +0+0=0.
\]
\end{proof}
\begin{remark}
Notice that the manifold $(P,\phi ,\theta , \xi)$ {\bf is not
cosymplectic} in sense of the definition of D.Blair (see \cite{BL},
Sec.6.5). More specifically, the (1,1)-tensor field $\phi $ is not
parallel.  To see this we notice that the formula
\[
2g((\nabla_{X}\phi)Y,Z)=g(N^{(1)}(Y,Z),\phi X)+2d\theta (\phi
Y,X)\theta (Z)-2d\theta (\phi Z,X)\theta (Y)
\]

for the covariant derivatives of the (1,1)-tensor $\phi$, proved in
\cite{BL}, Corollary 6.1, is true for an indefinite metric as well.
In order to see that $\phi $ is not parallel, substitute in this
formula $X=Z=X_{i},\ Y=\xi $.  Then $N^{1}=0$ by the previous
Proposition, the second term in the right side of this formula is
zero since $X_{i}$ is $\theta$-horizontal, and
\[
-2d\theta (\phi X_{i},X_{i})\theta (\xi )=-2\omega
(P_{i},X_{i})=-2.
\]
\par
Moreover, since $\theta$ is not a closed form, $(P,\phi ,\theta ,
\xi)$ is not cosymplectic in sense of P.Libermann either.
\end{remark}
\begin{remark} In the case of a contact metric manifold $(M,\theta
,g)$ one has $Ric(\xi )=2n-tr(h^2)$ for the tensor field
$h=\frac{1}{2}N_{J}^{(3)}$ with the component $N_{J}^{(3)}$ of the
Nijenhuis Tensor $N_{J}$, see \cite{BL}, Corollary 7.1.  Unlike
the Riemannian case, we have  for the metric $G$
$Ric(\xi)=R_{ij}\xi^{i}\xi^{j}=R_{00}=-\frac{n}{2}$.
\end{remark}

\vskip 1cm

\section{Group action of $ A_{1}^n$.}

Denote by $A_1$ the Lie group of affine transformations of the
real line $\mathbb{R}$.  In this section, we define an action of
the product $A_{1}^n$ on the space $\tP$ of symplectization which
is similar to the action of the Heisenberg group $H_{n}$ on $P$
(see \cite{OV}).\par

The Lie group $A_{1}$ of affine transformations of real line
$\mathbb{R}$ can be identified with the group of $2\times 2$ real
matrices of the form

\[
g=\begin{pmatrix} h & z\\ 0 & 1
\end{pmatrix}.
\]
The Lie algebra ${\mathfrak a}_{1}$ of the group $A_{1}$ in this
representation consists of matrices
\[
Y(a,z)=\begin{pmatrix} a & z\\ 0 & 0,
\end{pmatrix}
\]
and the exponential mapping $exp: {\mathfrak
\alpha}_{1}\rightarrow A_{1}$ takes the form
\[
exp\begin{pmatrix} a & z\\ 0 & 0
\end{pmatrix}=\begin{pmatrix} e^a & e^{a}z\\ 0 & 1
\end{pmatrix}.
\]
Left translations by the elements of the basic one-parameter group
$exp(tY(1,0))$:
\[
\begin{pmatrix} e^t & 0\\ 0 & 1
\end{pmatrix} \cdot \begin{pmatrix} h & z\\ 0 & 1
\end{pmatrix}=\begin{pmatrix} e^t h & e^t z\\ 0 & 1
\end{pmatrix}
\]
generate the basic right invariant vector field
\[
\xi_{a}=h\partial_{h}+z\partial_{z},
\]
while a similar action of the one-parameter group $exp(tY(0,1))$
produces the right invariant vector field
\[
\xi_{z}=\partial_{z}.
\]
We have $[\xi_{a},\xi_{z}]=-\xi_{z}$.

\par
For the right translations we have, respectively
\begin{multline}
\begin{pmatrix} h & z\\ 0 & 1
\end{pmatrix}\cdot \begin{pmatrix} e^t & 0\\ 0 & 1
\end{pmatrix}  =\begin{pmatrix} e^t h &  z\\ 0 & 1
\end{pmatrix}  \Rightarrow    \eta_{a}=h\partial_{h},\\
\begin{pmatrix} h & z\\ 0 & 1
\end{pmatrix}\cdot \begin{pmatrix} 1 & t\\ 0 & 1
\end{pmatrix}  =\begin{pmatrix}  h &  z+ht\\ 0 & 1
\end{pmatrix} \Rightarrow
\ \eta_{z}=h\partial_{z}.
\end{multline}

We have $[\eta_{a},\eta_{z}]=\eta_{z}$.\par

Consider now the identification of the Lie group $A_{1}$ with the
upper half-space $\mathbb{R}^2_{+}=\{(p,x)\in \mathbb{R}^2 \vert
p>0 \}$
\[
{\tilde \chi }: \begin{pmatrix}  h &  z\\ 0 & 1
\end{pmatrix} \Leftrightarrow \begin{pmatrix} p\\ x\end{pmatrix}=
\begin{pmatrix} h\\ -h^{-1}z\end{pmatrix} ,\ {\tilde \chi}^{-1}:
\begin{pmatrix} p\\
x\end{pmatrix}\Rightarrow \begin{pmatrix}  p &  -p^{-1}x\\ 0 & 1
\end{pmatrix}
\]
Under this identification
\[
\chi (\begin{pmatrix}  e^t &  0\\ 0 & 1
\end{pmatrix}\begin{pmatrix}  h &  z\\ 0 & 1
\end{pmatrix})=\begin{pmatrix} e^{t}h\\
-zh^{-1}\end{pmatrix},
\]
so that
\[
\chi_{*}(\xi_{a})=p\partial_{p}.
\]
Since
\[
\chi (\begin{pmatrix}  1 &  t\\ 0 & 1
\end{pmatrix}\begin{pmatrix}  h &  z\\ 0 & 1
\end{pmatrix})=\begin{pmatrix} h\\
-(z+t)h^{-1}\end{pmatrix},
\]
we get
\[
\chi_{*}(\xi_{c})=-p^{-1}\partial_{X}.
\]
We also calculate
\[
\chi^{*}(pdx)= -dz-zh^{-1}dh,\ \chi^{*}(dp\wedge
dx)=h^{-1}dh\wedge dz.
\]
\par
Now we apply these consideration to the product mapping
\[
\chi =\prod_{i=0}^{i=n}\chi_{i}: \prod_{i=0}^{i=n}A_{1}^{i}
\Leftrightarrow \tP_{+}=\{(p_{i},x_{i})\in \tP \vert p_{i}>0.\}
\]
We see that this mapping is the diffeomorphism satisfying

\[
\chi_{*}(\xi_{i\ a})=p_{i}\partial_{p_{i}}=P_{i},\
\chi_{*}(\xi_{i\ z})=-p^{-1}_{i}\partial_{x^{i}}=-L_{i}.
\]
In addition to this,
\begin{equation}
\begin{aligned}
\chi^{*}(\sum_{i}p_{i}dx^{i})&=-\sum_{i}(dz_{i}-z_{i}h_{i}^{-1}dh_{i}),\\
\chi^{*}(\sum_{i}dp_{i}\wedge
dx^{i})&=\sum_{i}p^{-1}_{i}dp_{i}\wedge dz_{i}.
\end{aligned}
\end{equation}
Dual to the right invariant frame $\xi_{a},\xi_{z}$ is the coframe
\[
\omega_{a}=h^{-1}dh,\ \omega_{z}=dz.
\]
Thus,
\[
\chi^{*}(\sum_{i}dp_{i}\wedge
dx^{i})=\sum_{i}p^{-1}_{i}dp_{i}\wedge dz_{i}
\]
is a right invariant symplectic structure on the group
$A_{1}^{n+1}.$\par
\begin{proposition}
The mapping
\[
\chi =\prod_{i=0}^{i=n}\chi_{i}: \prod_{i=0}^{i=n}A_{1}^{i}
\Leftrightarrow \tP_{+}=\{(p_{i},x_{i})\in \tP \vert p_{i}>0\}
\]
defined by
\[
\prod_{i=0}^{i=n}\begin{pmatrix}  h_{i} &  z_{i}\\ 0 & 1
\end{pmatrix} \Leftrightarrow \{  (p_{i}, x^{i})=
(h_{i},\  -h_{i}^{-1}z_{i})\}
\]
defines a symplectomorphism of the symplectic space $({\tilde
P}_{+}, \sum_{i=0}^{n}p_{i}dx^{i})$ with the product
$\prod_{i=0}^{i=n}A_{1}^{i}$ of (n+1) copies of affine group $A_{1}$
endowed with the symplectic structure generated by the right
invariant 1-form $-\sum_{i=0}^{i=n}(dz_{i}-z_{i}h_{i}^{-1}dh_{i}).$
\end{proposition}

\vskip 2cm

\section{Hyperbolic Rotations and the projectivization of $\tP$.}
In this section we construct a natural compactification of the TPS
$P$ endowed with the extension of the contact structure $\theta $
and that of the indefinite metric $G$.\par

Consider the action of the one-parameter group $R $ in the space
$\tP$ acting by the one-parameter group $HR$ of {\bf hyperbolic
rotations}
\begin{equation}
g^{t}: (p_{l},x^{i})\rightarrow (e^{t}p_{l},e^{-t}x^{i}).
\end{equation}
We have obviously
\begin{lemma}
\begin{enumerate}
\item The 1-form $\tilde \theta $ is invariant under this action
of the group HR. \item The metric $\tG$ is invariant under the
action of the group HR.
\end{enumerate}
\end{lemma}
\begin{proposition}
The space $\hat P$ of orbits of the points $\tP \setminus 0$ under
the action of the group HR is canonically isomorphic to the
projective space $P_{2n+1}(\mathbb{R})$.
\end{proposition}
\begin{proof}
Cover the space $\mathbb{R}^{2n+2}$ with the open subsets of two
types:
\par
Sets of the first type are
\[
{\hat U}_{j}=\{m\vert p_{j}\ne 0\},
\]
and associate with these sets the affine domains
\[
U_{j}\equiv \mathbb{R}^{2n+1}
\]
of the projective space $P_{2n+1}(\mathbb{R})$ with the
coordinates
\[
(x^{i}p_{j},\frac{p_{l}}{p_{j}}).
\]
Sets of the second type are
\[
{\hat V}_{k}=\{m\vert x^{k}\ne 0\},
\]
and associate with these sets the affine domains
\[
V_{k}\equiv \mathbb{R}^{2n+1}
\]
of the projective space $P_{2n+1}(\mathbb{R})$ with the
coordinates
\[
(\frac{x^{i}}{x^{k}},p_{l}x^{k}).
\]
Notice that on the intersections $U_{j_{1}}\cap U_{j_{2}}$ we have
relations between the corresponding affine coordinates
\[
x^{i}p_{j_{2}}=x^{i}p_{j_{1}}\cdot (\frac{p_{j_{2}}}{p_{j_{1}}});\
\ \frac{p_{k}}{p_{j_{2}}}=\frac{p_{k}}{p_{j_{1}}}\cdot
(\frac{p_{j_{1}}}{p_{j_{2}}}).
\]
On the intersections $U_{j}\cap V_{k}$, we have relations between
the corresponding affine coordinates
\[
x^{l}p_{j}=\frac{x^{l}}{x^{k}}\cdot (x^{l}p_{j});\ \
\frac{p_{l}}{p_{j}}=p_{l}x^{k}\cdot (\frac{1}{x^{k}{p_{j}}}).
\]
Finally, on the intersections $V_{j_{1}}\cap V_{j_{2}}$ we have
relations between the corresponding affine coordinates
\[
\frac{x^{k}}{x^{j_{2}}}=\frac{x^{k}}{x^{j_{1}}}\cdot
(\frac{x^{j_{1}}}{x^{j_{2}}});\ \
p_{l}x^{j_{2}}=p_{l}x^{j_{1}}\cdot (\frac{x^{j_{2}}}{x^{j_{1}}}).
\]
This shows that affine coordinates of all the affine charts are
related by the transition functions {\bf invariant under the
action of hyperbolic rotations}.   Thus, they are glued into the
standard projective space $P_{2n+1}(\mathbb{R})$.
\end{proof}
Combining Lemma 2 and the previous construction we get the
following
\begin{proposition}
\begin{enumerate}
\item The projections $\hat \theta$ of the 1-form $\tilde \theta$  and
that of the metric ${\hat G}$ of $\tG$ endow the projective space
${\hat P}$ with the contact structure and the metric of signature
$(n+1,n).$
\item The composition $J$ of the embedding $j:P\rightarrow \tP$
and the projection $\tP \rightarrow {\hat P}$ defined the {\bf
compactification} $({\hat P},{\hat \theta},{\hat G})$ of the TPS
$(P,\theta ,G)$ with the contact structure and Mrugala metric $G$.
\end{enumerate}
\end{proposition}
\section{Group action of $H_{n}$ and the "partial orbit structure" of $\hat
P$.}

In this section we consider the lift to the space $\tP$ of the
action of the group $H_{n}$ on $P$ discussed in Sec. 12 and the
action of subgroups of $H_{n}$ on the cells of smaller dimension
of the standard CW-structure of the projective space $\hat P$.\par

The differential operators $X_{i},P_{j}, \xi $ of the canonical
frame (5.5) act also in the space $\tP$ with the same commutator
relations. This action generates the action of the Lie group
$H_{n}$ on the space $\tP$ leaving hyperplanes $p_{0}=const $
invariant.
\par

Introduce the sequence $L_{k},\ k=1,\ldots ,n$ of subgroups of the
Heisenberg group $H_{n}$ defined by the condition
\begin{equation}
L_{n-k}=\{ g({\bar a},{\bar b},c)\vert b_{1}=\ldots = b_{k}=0\}.
\end{equation}
These subgroups form the series
\[
H_n \supset L_{n-1}\supset L_{n-2}\supset \ldots \supset L_{0}.
\]
It is easy to see that
\begin{equation}
L_{n-k}\simeq \mathbb{R}^{k}\times H_{n-k}
\end{equation}
is the product of the $k$-dim abelian group $\mathbb{R}^k$ and the
Heisenberg group $H_{n-k}$. \par

The right invariant vector fields on $H_{n}$ tangent to (and
generated by) the subgroup $L_{n-k}$ are (in terms of the
isomorphism of Sec. 12) $\xi , X_{i},i=1,\ldots ,n, \
P_{j},j=k+1,\ldots ,n$.
\par

In the space $\tP$, consider the affine planes $V_{k}$ defining
the cells of the standard cell structure of the projective space
${\hat P}=P_{2n+1}(\mathbb{R})$ with respect to the
(hyperbolically) homogeneous coordinates of $\hat P$
\begin{equation}
V_{k}=\{ (x^{i},p_{j}\vert p_{0}=p_{1},\ldots =p_{k-1}=0,\
p_{k}=1\},
\end{equation}
with $k=0,1,\ldots n$.  It is clear that the projective space
$\hat P$ is obtained by gluing to the cell $V_{0}=j(P)$ the
smaller cells $V_{1},V_2,V_3,$ consecutively and, finally, by
gluing in n-dim projective space $P_{n}(\mathbb{R})$ obtained by
the action of hyperbolic rotations (usual dilatations here) on the
subspace $V_{n+1}=\{(x^{i},p_{j}=0\vert j=0, \ldots ,n\}.$\par

It is easy to see now that each cell $V_{k}$ is {\bf canonically
diffeomorphic} to the group $L_{n-k}$ whose action on $V_{k}$ is
induced by the action of the Heisenberg group $H_{n}$ on the space
$\tP$ considered above.
\par

So, even though the action of $H_{n}$ on $P$ cannot be extended to
the compactification $\hat P$, a coherent action of the subgroups
of series (21.1) produces the partial cell structure of $\hat P$
starting with the projective subspace $V_{n+1}\simeq
P_{n}(\mathbb{R})$.
\par

The restriction of the 1-form $\tilde \theta$ to the cell $V_{k}$
has the form
\[
\theta_{k}={\tilde
\theta}\vert_{V_{k}}=dx^{k}+\sum_{i=k+1}^{n}p_{i}dx^{i}.
\]
Therefore, this form determines the canonical contact structure on
the Heisenberg factor of the cell $V_{k}\simeq L_{n-k}\simeq
\mathbb{R}^{k}\times H_{n-k}$ and is zero on the first factor.\par

The restriction of the Mrugala metric $G$ to the cell $V_{k}$ has,
in variables $(p_{k},\ldots ,p_{n};x^{0},\ldots,
x^{k-1};x^{k},\ldots , x^{n})$, the form
\begin{equation}
G_{k}={\tG}\vert_{V_{k}}=\begin{pmatrix} 0_{(n-k)\times (n-k)} &
0_{(n-k)\times k} & I_{n-k}\\
0_{k\times (n-k)} & 0_{k \times k} & 0_{k\times (n-k)}\\
I_{n-k} & 0_{(n-k)\times (n-k)} & p_{i}p_{j}
\end{pmatrix}
\end{equation}
Thus, this metric is zero on the first abelian factor of the cell
$V_{k}\simeq L_{n-k}$ and coincides with the Mrugala metric on the
Heisenberg factor $H_{n-k}$ of the cell $V_{k}$.\par Combining these
arguments we get the following
\begin{theorem}
\begin{enumerate}
\item The restriction of the action of the Heisenberg group
$H_{n}$ on the space $P\simeq V_{0}$ (embedded in $\hat P$) to the
subgroup $L_{n-k}\simeq R^{k}\times H_{n-k}$ of the form (23.1)
extends to the action of this subgroup on the cell $\pi
(V_{k})\subset {\hat P}$ and determines the diffeomorphism of the
group $L_{n-k}$ with $V_{k}$ and with its image $\pi (V_{k})
\subset {\hat P}.$

\item Restrictions of the 1-form $\tilde \theta $ and metric $\tG$
to the cell $V_{k}$ endow the Heisenberg factor $H_{n-k}$ of
$V_{k}$ with the contact structure $\theta_{k}$ and the Mrugala
metric $G_{k}$ and are both zero on the abelian factor
$\mathbb{R}^k$.
\end{enumerate}
\end{theorem}

\begin{remark}
Every cell $V_{k}$ represents the thermodynamical phase space of
an abstract thermodynamical system with $n$ extensive and $n-k$
intensive variables. This corresponds to a situation where the
thermodynamical potential $\phi (x^{i})$ depends on
$x_{k+1},\ldots x^{n}$ but not on the first $k$ extensive
variables $x^i$.  As a results the participation of factors $x^i
,i=0,\ldots, k-1$ in the processes is "switched out" and they
become parameters only.\par
\end{remark}
\begin{example}
Consider the case $n=2$, i.e. take $P^5$ to be five-dimensional
with the contact form $\theta =dU-SdT+pdV$ (a one-component
homogeneous system, per 1 mole). Hypersurface $\mathcal{C}$ (see
Sec.11) has, in this case, the well known form $U-ST+pV=0$.  Its
lift to $\tP$ - ${\tilde {\mathcal{C}}}$ has the form
$p_{0}U-ST+pV=0.$\par

 Intersection of this quadric with the plane
$p_{0}=0$ is the (degenerate) quadric $\mathcal{C}_{1}:\ pV=ST$.
Fixing value of $S$, say, taking $S=R-const$ determine the cell
$\simeq V_{1}$ that projects onto the cell $V_{1}$ of the compact
space $\hat P$. Image of the quadric $\mathcal{C}_{1}$ under this
projection determine in the 3dim $H_{1}$-factor of the cell
$V_{1}$ the surface $pV=RT$ given by the equation of mono-atomic
ideal gas.  \par

Hypersurface $\mathcal{C}$ is the submanifold containing all the
constitutive (equilibrium) surfaces of all thermodynamical systems
with the TPS $P^5$.  Closures of these surfaces in  $\hat P$
contains points from cells of smaller dimension $V_{1},V_{2}.$
Thus, equation of a mono-atomic ideal gas appears here as the
equation of the surface formed by the limit points in $V_{1}$ of
all possible constitutive surfaces in $(P^5,\theta =dU-SdT+pdV)$.
\end{example}
\begin{remark}
The construction of a compact manifold $\hat P$ in terms of a
series of subgroups (23.1) of a Lie group $H$ represents a way to
represent a manifold in terms of a Lie group with the open dense
orbit isomorphic to the group $H$ itself and a natural Whitney
stratification in terms of extension of subgroup actions on the
cells of smaller dimension.  Removed generators (here $P_{i},\
i=1,\ldots ,k$) determine the projections from cells of higher
dimension to the cells in their closure (see \cite{T}).
\end{remark}

\vskip 1cm
\section{Conclusion}
In this work we've studied basic properties of the indefinite metric
$G$ of R.Mrugala defined on the contact (2n+1)-dimensional phase
space $(P,\theta )$ of a homogeneous thermodynamical system. We have
calculated the curvature tensor, Killing vector fields, and the
second fundamental form of the Legendre submanifolds of $P$ -
constitutive surfaces of different homogeneous thermodynamical
systems. We established an isomorphism of the TPS $(P,\theta ,G)$
with the Heisenberg Lie group $H_{n}$ endowed with the right
invariant contact structure and the right invariant indefinite
metric. We lifted the metric $G$ to a metric $\tG$ of signature
$(n+1,n+1)$ in the symplectization $\tP$ of the contact space
$(P,\theta)$ and studied curvature properties and Killing vector
fields of this metric. Finally we introduced the "hyperbolic
projectivization" of the space $(\tP,{\tilde \theta}, \tG)$ that can
be considered as the natural {\bf compactification} (with the
contact structure and the indefinite metric) of the TPS space
$(P,\theta ,G).$  \par

Many interesting questions were left outside of this paper - study
of the geodesics of metric $G$, the relation of metric properties
of $(P,G)$ with the contact transformations (see \cite{MNSS, Be}),
the characterization of the submanifolds of signature changes of
the thermodynamical metrics on the Legendre submanifolds (related
to the phase transitions in the corresponding homogeneous
thermodynamical systems) in terms of the Grassmanian of the
Legendre n-dim subspaces in $P$, the use of the geometry of
Heisenberg group $H_{n}$ to the study of Legendre submanifolds of
$(P,\theta)$, and others.  Some of these questions will be
considered in the continuation of this work.\par

In the conclusion we would like to thank Professor M. Goze for the
useful information about contact and metric structures on the
Heisenberg group.

\section{Appendix: Killing vector fields for $\tG $.}
In this Appendix we provide the details of calculations of Killing
vector fields of metric $\tG $.\par

Conditions for a vector field $Z$ on a manifold $M$ with a
(pseudo-Riemannian) metric $g$ to be a Killing vector field have
the form ${\mathcal L}_{Z}g=0$, or, in local coordinates
$(p_{i},x^{j})$,
\begin{equation}
\frac{\partial g_{ij}}{\partial z^{k}}Z^{k}+g_{ik}\frac{\partial
Z^{k}}{\partial z^{j}}+g_{kj}\frac{\partial Z^{k}}{\partial
z^{i}}=0,\ \text{for all}\ i,j,k.
\end{equation}
\par We
take a vector field $X$ in $\tP$ in the form
\begin{equation}
X =\xi^{p_{l}}\partial_{p_{l}}+\xi^{x^{s}}\partial_{x^{s}}
\end{equation}
and consider cases of different pairs of indices $(ij)$.\par
Case $(ij)=(p_{l}p_{s})$.\par
\[
0
+\delta^{p_{l}}_{k}\partial_{p_{s}}\xi^{k}+\delta^{p_{s}}_{k}\partial_{p_{l}}\xi^{k}=0
\]
or
\begin{equation}
\partial_{p_{s}}\xi^{x^{l}}+\partial_{p_{l}}\xi^{x^{s}}=0.
\end{equation}
Case $(ij)=(p_{l}x^{s}).$
\[
0+\delta^{l}_{k}\partial_{x^{s}}\xi^{x^{k}}+
(\delta^{s}_{k}\partial_{p_{l}}\xi^{p_{k}}+p_{s}p_{k}\partial_{p_{l}}\xi^{x^{k}})=0,
\]
or
\begin{equation}
\partial_{x^{s}}\xi^{x^{l}}+
\partial_{p_{l}}\xi^{p_{s}}+p_{s}p_{k}\partial_{p_{l}}\xi^{x^{k}}=0.
\end{equation}
Case $(ij)=(x^{i}x^{j})$.
\[
(p_{i}\xi^{p_{j}}+p_{j}\xi^{p_{i}})+
(\delta^{i}_{k}\partial_{x^{j}}\xi^{p_{k}}+p_{i}p_{k}\partial_{x^{j}}\xi^{x^{k}})+
(\delta^{j}_{k}\partial_{x^{i}}\xi^{p_{k}}+p_{j}p_{k}\partial_{x^{i}}\xi^{x^{k}})=0,
\]
or
\begin{equation}
(p_{i}\xi^{p_{j}}+p_{j}\xi^{p_{i}})+
(\partial_{x^{j}}\xi^{p_{i}}+p_{i}p_{k}\partial_{x^{j}}\xi^{x^{k}})+
(\partial_{x^{i}}\xi^{p_{j}}+p_{j}p_{k}\partial_{x^{i}}\xi^{x^{k}})=0.
\end{equation}
We will need the following
\begin{lemma}
Let $f^{i}(y^{1},\ldots y^{n})$ be $n$ functions of $n$ variables
$x^{i}$ such that for all $i,j$
\[
\partial_{y^{i}}f^{j}+\partial_{y^{j}}f^{i}=0,
\]
then
\[
f^{i}=a^{i}_{j}y^{j}+b^{i}
\]
with some constant vector $b^{i}$ and skew-symmetric constant
matrix $(a^{i}_{j})$ : $a^{j}_{i}=-a^{i}_{j}$.
\end{lemma}
\begin{proof}
Apply $\partial_{y^k}$ to the condition of Lemma.  We get
\[
\partial_{y^{k}}\partial_{y^{i}}f^{j}+\partial_{y^{k}}\partial_{y^{j}}f^{i}=0.
\]
cyclic permutation of indices $kij$ give us two more equalities.
Add first two equations and subtract the third one.  We get
\[
2\partial_{y^{i}}\partial_{y^{j}}f^{k}=0.
\]
This being true for all triples of indices $ijk$ shows that {\bf
all} second derivatives of all functions $f^k$ are zero.\par
Therefore
\[
f^{i}=\sum_{j}a^{i}_{j}y^{j}+b^{i}
\]
with constant coefficients.\par
Writing down condition of Lemma for
these linear functions we find that the matrix $a^{i}_{j}$ is
skew-symmetric.
\end{proof}

Applying this Lemma to the functions $\xi^{x^{i}}$ with $p_{l}$ as
arguments in Lemma we have due to the equality (23.3)
\begin{equation}
\xi^{x^{i}}=a_{i}^{l}(x)p_{l}+b^{i}(x)
\end{equation}
with an skew-symmetric matrix function  $a^{i}_{j}(x)$ and scalar
functions $b^{i}(x)$.  This solves equations (23.3).\par

Substituting these expressions into the other two families of
equations we present these equations in the form
\begin{multline}
\partial_{x^{s}}a^{k}_{l}(x)p_{k}+\partial_{x^{s}}b^{l}+\partial_{p_{l}}\xi^{p_{s}}+p_{s}p_{k}a^{l}_{k}(x)=0,\\
p_{i}\xi^{p_{j}}+p_{j}\xi^{p_{i}}+(\partial_{x^{j}}\xi^{p_{i}}+\partial_{x^{i}}\xi^{p_{j}})+
(p_{i}p_{k}[\partial_{x^{j}}(a_{k}^{l}(x)p_{l}+b^{k}(x))]+
p_{j}p_{k}[\partial_{x^{i}}(a_{k}^{l}(x)p_{l}+b^{k}(x))])=0.
\end{multline}
Rewrite first equation in the form
\begin{equation}
\partial_{p_{l}}\xi^{p_{s}}=-\partial_{x^{s}}a^{k}_{l}(x)p_{k}+p_{s}p_{k}a^{k}_{l}(x)-\partial_{x^{s}}b^{l},
\end{equation}
where we have used skew-symmetry of matrix $a$ in the last term,
and apply $\partial_{p_{m}}$ to this formula. We get
\[
\partial_{p_{m}}\partial_{p_{l}}\xi^{p_{s}}=-\partial_{x^{s}}a^{m}_{l}(x)+p_{s}a^{m}_{l}(x)+\delta_{sm}a^{k}_{l}p_{k}.
\]
Switching $m$ and $l$ we get
\[
\partial_{p_{l}}\partial_{p_{m}}\xi^{p_{s}}=-\partial_{x^{s}}a^{l}_{m}(x)+p_{s}a^{l}_{m}(x)+\delta_{sl}a^{k}_{m}p_{k}.
\]
Equating terms in the right side that do not depend on $p$ we get
\[
\partial_{x^{s}}a^{m}_{l}(x)=\partial_{x^{s}}a^{l}_{m}(x)\rightarrow
\partial_{x^{s}}(2a^{m}_{l}(x)))=0.
\]
Thus, matrix $(a^{m}_{l})$ is constant.\par Equating terms linear
by $p$ we get
\[
p_{s}a^{m}_{l}(x)+\delta_{sm}a^{k}_{l}p_{k}=p_{s}a^{l}_{m}(x)+\delta_{sl}a^{k}_{m}p_{k}.
\]
For $s\ne m,l$ this gives $a^{l}_{m}=a^{m}_{l}$ and, due to the
skew-symmetry of matrix $a$, $a^{l}_{m}=0$.  \par For $s=l\ne m$
we get $p_{l}a^{m}_{l}=p_{l}a^{l}_{m}+a_{m}^{k}p_{k}$.  From this
it follows that $a_{m}^{k}=0$ for all $k\ne l$ and that
$a_{m}^{l}=0$ as well. Finally, if $s=l=m$ we get by the
skew-symmetry of $a^{i}_{j}$ that $a^{l}_{l}=0$. Therefore, matrix
$a$ is zero $(a^{i}_{j})=0$ and
\begin{equation}
\xi^{x^{s}}=b^{s}(x).
\end{equation}
\par
As a result, first of the equations (23.7) reads now
\[
\partial_{p_{l}}\xi^{p_{s}}=-\partial_{x^{s}}b^{l}(x),
\]
i.e. $\xi^{p_{s}}$ depends on variables $p_{l}$ linearly
\begin{equation}
\xi^{p_{s}}=-\partial_{x^{s}}b^{l}(x)p_{l}+h^{s}(x).
\end{equation}
This solves the first of equations (23.7).
\par

The second equation in (23.7) can now be written in the form
\begin{multline}
p_{i}(-\partial_{x^{j}}b^{l}(x)p_{l}+h^{j}(x))+p_{j}(-\partial_{x^{i}}b^{l}(x)p_{l}+h^{i}(x))+
(\partial_{x^{j}}(-\partial_{x^{i}}b^{l}(x)p_{l}+h^{i}(x))+\\
+\partial_{x^{i}}(-\partial_{x^{j}}b^{l}(x)p_{l}+h^{j}(x)))+
(p_{i}\partial_{x^{j}}p_{k}b^{k}(x)+
p_{j}\partial_{x^{i}}p_{k}b^{k}(x))=0,
\end{multline}
or
\begin{equation}
p_{i}h^{j}(x)+p_{j}h^{i}(x)+(\partial_{x^j}h^{i}(x)+\partial_{x^{i}}h^{j}(x))-
  2\partial_{x^{i}}\partial_{x^{j}}(b^{l}(x)p_{l})=0.
\end{equation}

Terms independent on $p$ give us
\[
\partial_{x^j}h^{i}(x)+\partial_{x^{i}}h^{j}(x)=0,
\]
for all $i,j$.  The lemma above gives us
\begin{equation}
h^{i}(x)=\sum_{j}q^{i}_{j}x^{j}+k^{i}
\end{equation}
with constant coefficients and $q^{i}_{j}=-q^{j}_{i}.$\par

Linear by $p$ part of equality (23.12) has the form
\begin{equation}
p_{i}h^{j}(x)+p_{j}h^{i}(x)=
  2\partial_{x^{i}}\partial_{x^{j}}(b^{l}(x)p_{l}).
\end{equation}
Taking derivatives by $p_{i}$ and then by $x^{l}$ we get,

\begin{equation}
\partial_{p_{i}}:\hskip1cm
h^{j}(x)+\delta^{j}_{i}h^{i}(x)=2\partial_{x^{i}}\partial_{x^{j}}b^{i}(x)\hskip
0.5cm \text{no summation},
\end{equation}
and then as the coefficient of $p_{i}$ (using (23.13):
\begin{equation}
\partial_{x^{l}}:\hskip1cm
q^{j}_{l}+\delta^{j}_{i}q^{i}_{l}=2\partial_{x^{l}}\partial_{x^{i}}\partial_{x^{j}}b^{i}(x)\hskip
0.5cm \text{no summation}.
\end{equation}
Since the right side in the second equation(s) is symmetrical by
$(l,j)$, one has
\[
q^{j}_{l}+\delta^{j}_{i}q^{i}_{l}=q^{l}_{j}+\delta^{l}_{i}q^{i}_{j},\
\hskip 1cm \text{no summation}.
\]
For $i=j$ last equality gives (since matrix $q$ is skew-symmetric)
\[
2q^{i}_{l}=q^{l}_{i}=-q^{i}_{l},
\]
so that for all $i,l$ $q^{i}_{l}=0$ and therefore,
\begin{equation}
h^{i}(x)=k^{i}.
\end{equation}
Substituting this back to (23.14) we get
\begin{equation}
p_{i}k^{j}+p_{j}k^{i}=
  2\partial_{x^{i}}\partial_{x^{j}}(b^{l}(x)p_{l}).
\end{equation}
Differentiating by $p_{l},\ l\ne i,j$ we get
\[
\partial_{x^{i}}\partial_{x^{j}}b^{l}(x)=0.
\]
As a result,
$b^{l}(x)$ is linear by all variables except $x^{l}$
\begin{equation}
b^{l}(x)=\sum_{k\ne l}\alpha^{l}_{k}(x^{l})x^{k}+\beta^{l}(x^{l}).
\end{equation}
Using $h^{i}(x)=k^i$ in (23.15) we get
\begin{equation}
k^{j}+\delta^{j}_{i}k^{i}=2\partial_{x^{i}}\partial_{x^{j}}b^{i}(x)\hskip
0.5cm \text{no summation}.
\end{equation}
Taking here first $i=j=l$, then $i=l\ne j=r$ we find
\[
\partial^{2}_{x^{l}}b^{l}(x)=k^l ,\
\partial_{x^{r}}\partial_{x^{l}}b^{l}(x)=\frac{1}{2}k^{r}.
\]
Substituting here the expression for $b^{l}(x)$ from (23.9) we get
the first equation in the form
\[
k^{l}=\sum_{r\ne
l}\partial^{2}_{x^{l}}\alpha^{l}_{r}x^{r}+\partial^{2}_{x^{l}}\beta^{l}(x^{l}).
\]
Therefore, $\alpha^{l}_{r}$ is {\bf linear} by $x^{l}$ while
\[
\beta^{l}=\frac{1}{2}k^{l}(x^{l})^{2}+\mu^{l}x^{l}+\nu^{l}.
\]
The second equation, valid for $l\ne r$, gives
$\frac{1}{2}k^{l}=\partial_{x^{l}}\alpha^{l}_{r}(x^{l})$ and
therefore
\[
\alpha^{l}_{r}=\frac{1}{2}k^{r}x^{l}+\eta^{l}_{r}.
\]
Combining yields the coefficients of $b^{i}$ and using them in
(23.19) we find
\begin{equation}
\xi^{x^{l}}=\sum_{r\ne
l}(\frac{1}{2}k^{r}x^{l}+\eta^{l}_{r})x^{r}+\frac{1}{2}k^{l}(x^{l})^{2}+\mu^{l}x^{l}+\nu^{l}=
\sum_{r}(\frac{1}{2}k^{r}x^{r})x^{l}+\sum_{r\ne
l}\eta^{l}_{r}x^{r}+\mu^{l}x^{l}+\nu^{l},
\end{equation}
or, renaming coefficient $\mu^{l}=\eta^{l}_{l}$,
\begin{equation}
\xi^{x^{l}}=\sum_{r}(\frac{1}{2}k^{r}x^{r})x^{l}+\sum_{r}\eta^{l}_{r}x^{r}+\nu^{l}.
\end{equation}
The equality (23.10) with the found values of the coefficients gives
us
\begin{multline}
\xi^{p_{s}}=-\sum_{l}\partial_{x^{s}}\xi^{x^{l}}p_{l}+k_{s}=k_{s}-\sum_{l}(\eta^{l}_{s}+
\frac{1}{2}\delta^{l}_{s}\sum_{r}(k^{r}x^{r})+\frac{1}{2}x^{l}k^{s})p_{l}=\\
-\sum_{l}\eta^{l}_{s}p_{l}+\sum_{r}k^{r}\left(\delta_{r}^{s}-
\frac{1}{2}\delta_{r}^{s}(\sum_{l}x^{l}p_{l})-\frac{1}{2}x^{r}p_{s}
\right).
\end{multline}
As a result we get for a Killing vector field $X$ the following
representation with arbitrary scalar coefficients $k^r ,
\eta^{s}_{r} , \nu^{s}$
\begin{multline}
X =\xi^{x^{s}}\partial_{x^{s}}+\xi^{p_{l}}\partial_{p_{l}}=
\sum_{s} \left(
\sum_{r}(\frac{1}{2}k^{r}x^{r})x^{s}+\sum_{r}\eta^{s}_{r}x^{r}+\nu^{s}\right)
\partial_{x^{s}}+\\ +\sum_{s} \left(
-\sum_{l}\eta^{l}_{s}p_{l}+\sum_{r}k^{r}\left(\delta_{r}^{s}-
\frac{1}{2}\delta_{r}^{s}(\sum_{l}x^{l}p_{l})-\frac{1}{2}x^{r}p_{s}
\right)\right)\partial_{p_{s}}.
\end{multline}

Splitting these expressions in accordance with the different
independent parameters we get the following basis of the Lie algebra
of Killing vector fields:

$\nu^s$-term
\begin{equation}
X_{s}=\partial_{x^{s}},
\end{equation}

$\eta^{s}_{r} $ - term
\begin{equation}
Q^{r}_{s}=x^{r}\partial_{x^{s}}-p_{s}\partial_{p_{r}},
\end{equation}

$k^{i}$-terms
\begin{equation}
D^{i}=\frac{x^{i}}{2}\left[
x^{s}\partial_{x^{s}}-p_{l}\partial_{p_{l}}\right]
+(1-\frac{1}{2}(x^{l}p_{l}))\partial_{p_{i}}=
\frac{x^{i}}{2}Q+(1-\frac{1}{2}(x^{l}p_{l}))\partial_{p_{i}},
\end{equation}

where
\[Q=\sum_{s}Q^{s}_{s}=\sum_{s}(x^{s}\partial_{x^{s}}-p_{s}\partial_{p_{s}}).\]\par

Vector fields $X^i$ form the abelian Lie subalgebra
$\mathfrak x $ of Lie algebra of Killing vector fields ${\mathfrak
k}_{\tG}$.\par

It is easy to see that $(n+1)^2$ vector fields $Q^{r}_{s}$ form
the Lie subalgebra ${\mathfrak q}$ of the type ${\mathfrak
gl}(n+1,\mathbb{R})$ since
\[
[Q^{i}_{j},Q^{p}_{k}]=\delta^{p}_{j}Q^{i}_{k}-\delta^{i}_{k}Q^{p}_{j}
\]
similar to the commutator relations of the basic matrices
$E^{i}_{j}$ of ${\mathfrak gl}(n+1,\mathbb{R})$.  Vector field $Q$
generate the center of this Lie subalgebra.\par

We also have
\[
[Q^{i}_{j},X_{s}]=-\delta^{i}_{s}X_{j}, [Q,X_{s}]=-X_{s},
\]
so that adjoint action of ${\mathfrak q}$ on $\mathfrak x $ is
isomorphic to the  standard action of ${\mathfrak
gl}(n+1,\mathbb{R})$ on $\mathbb{R}^{n+1}$.\par

Calculate commutator $[D^{i},D^{j}]$
\begin{multline}
[D^{i},D^{j}]=[\frac{x^{i}}{2}Q
+(1-\frac{1}{2}(x^{l}p_{l}))\partial_{p_{i}}, \frac{x^{j}}{2}Q
+(1-\frac{1}{2}(x^{l}p_{l}))\partial_{p_{j}}]=\frac{x^{i}}{2}\{
\sum_{s}x^{s}(\delta^{j}_{s}Q+\frac{x^{j}}{2}\partial_{x^{s}})+\sum_{l}p_{l}\frac{x^{j}}{2}\partial_{p_{l}}+\\
+\sum_{s}x^{s}(-\frac{p_{s}}{2})\partial_{p_{j}}-\sum_{m}p_{m}(-\frac{x^{m}}{2})\partial_{p_{j}}\}
+(1-\frac{{\bar x}\cdot {\bar p}}{2})\{
\frac{x^{j}}{2}(-\partial_{p_{i}})-\frac{x^{i}}{2}\partial_{p_{j}}\}=\\
\frac{x^{i}}{2}\{\frac{1}{2}Qx^{j}+\frac{x^{j}}{2}\sum_{s}x^{s}\partial_{x^{s}}+\frac{x^{j}}{2}\sum_{l}p_{l}\partial_{p_{l}}
\}+ (1-\frac{{\bar x}\cdot {\bar
p}}{2})\{-\frac{1}{2}(x^{j}\partial_{p_{i}}+x^{i}\partial_{p_{j}})\}=\\
=\frac{1}{4}Qx^{i}x^{j}+\frac{1}{4}x^{i}x^{j}\sum_{s}x^{s}\partial_{x^{s}}+\frac{x^{i}x^{j}}{2}\sum_{l}p_{l}\partial_{p_{l}}+
(1-\frac{{\bar x}\cdot {\bar
p}}{2})\{-\frac{1}{2}(x^{j}\partial_{p_{i}}+x^{i}\partial_{p_{j}})\}=0.
\end{multline}
Last expression is equal zero since all its terms are {\bf
symmetrical} by $(ij)$.  To check this one can also use the
relations

\[
Q(1-\frac{<p,x>}{2})=0,\ [x^{i}Q,x^{j}Q]=0.
\]

\par

Thus, vector fields $D^{i}$ also form the abelian Lie subalgebra
$\mathfrak d$ of $\mathfrak k$.\par Next we calculate
\begin{multline}
[D^{s},Q^{i}_{j}]=[\frac{x^{s}}{2}Q
+(1-\frac{1}{2}(x^{l}p_{l}))\partial_{p_{s}},x^{i}\partial_{x^{j}}-p_{j}\partial_{p_{i}}]=
\frac{x^{s}}{2}\{
\sum_{m}x^{m}\delta^{m}_{i}\partial_{x^{j}}+\sum_{l}p_{l}\delta^{j}_{l}\partial_{p_{i}}\}+\\
+\{ (1-\frac{1}{2}(x^{l}p_{l}))\delta^{j}_{s}(-\partial_{p_{i}})\}
-
x^{i}[\delta^{s}_{j}\frac{1}{2}Q+\frac{x^{s}}{2}\partial_{x^{j}}
+(-\frac{1}{2}p_{j})\partial_{p_{s}}]+p_{j}[\frac{x^{s}}{2}(-\partial_{p_{i}})-\frac{x^{i}}{2}\partial_{p_{s}}]=\\
\frac{x^{s}}{2}\{
x^{i}\partial_{x^{j}}+p_{j}\partial_{p_{i}}\}-(1-\frac{{\bar
x}\cdot {\bar p}}{2})\delta^{j}_{s}\partial_{p_{i}}\}-
x^{i}[\delta^{s}_{j}\frac{1}{2}Q+\frac{x^{s}}{2}\partial_{x^{j}}
-\frac{p_{j}}{2}\partial_{p_{s}}]-\frac{x^{s}p_{j}}{2}\partial_{p_{i}}-\frac{x^{i}p_{j}}{2}\partial_{p_{s}}=\\
=-(1-\frac{{\bar x}\cdot {\bar
p}}{2})\delta^{j}_{s}\partial_{p_{i}}-\delta^{s}_{j}\frac{x^{i}}{2}Q=-\delta^{s}_{j}D^{i},
\end{multline}
or
\begin{equation}
[D^{s},Q^{i}_{j}]=-\delta^{s}_{j}D^{i},
\end{equation}
so, abelian subalgebra $\mathfrak d$ is invariant under the
adjoint action of $\mathfrak q$ with the same standard action as
for $\mathfrak x$.\par

Notice that
\[
Q\vert_{\mathfrak d}=Id_{\mathfrak d},\ Q\vert_{\mathfrak
x}=-Id_{\mathfrak x}.
\]
\par

Finally,

\begin{multline}
[X^{s},D^{i}]=[\partial_{x^{s}},\frac{x^{i}}{2}Q
+(1-\frac{1}{2}(x^{l}p_{l}))\partial_{p_{i}}]=
\frac{1}{2}\delta^{i}_{s}(\sum_{k}x^{k}\partial_{x^{k}}-\sum_{l}p_{l}\partial_{p_{l}})
+\frac{x^{i}}{2}\partial_{x^{s}}-\frac{1}{2}p_{s}\partial_{p_{i}}=\\
=\frac{1}{2}(x^{i}\partial_{x^{s}}-p_{s}\partial_{p_{i}})+\frac{1}{2}\delta^{i}_{s}Q=\frac{1}{2}Q^{i}_{s}
+\frac{1}{2}\delta^{i}_{s}Q.
\end{multline}
\par
Next step is to prove that the Lie algebra ${\mathfrak iso}_{\tG}$
of killing vector fields is isomorphic to the Lie algebra
${\mathfrak sl}(n+2,\mathbb{R}).$\par

Consider the Lie algebra ${\mathfrak gl}(n+2,\mathbb{R})$ of real
matrices and introduce several subspaces of this Lie algebra:\par

The Lie subalgebra $g_{n+1}$ of matrices
\[
g=\begin{pmatrix} A & 0\\
0 & 0
\end{pmatrix} , \ A\in {\mathfrak gl}(n+1,R)
\]
with the basis formed by matrices $E^{i}_{j}$ having the entry $1$
at the $ij$ place and zero otherwise,\par

the subspace $d$ of matrices of the form
\[
\begin{pmatrix}
0 & {\bar f}\\
0 & 0
\end{pmatrix},\ {\bar f}\in \mathbb{R}^n
\]
with the basis $\{ 1^{i} \}$ of matrices having entry $1$ at the
place $i(n+2)$ and zero at all other places,\par

subspace $x$ of matrices of the form
\[
\begin{pmatrix}
0 & 0\\
{\bar e} & 0
\end{pmatrix},\ {\bar e}\in \mathbb{R}^n
\]
with the basis $\{ 1_{j} \}$ of matrices having entry $1$ at the
place $(n+2)j$ and zero at all other places,\par

one-dimensional subspace of matrices
\[
\begin{pmatrix}
0 & 0\\
0 & t
\end{pmatrix},\ t\in \mathbb{R}
\]
with the basis $E_{**}$ formed by the matrix
$E_{**}=\begin{pmatrix}
0 & 0\\
0 & 1
\end{pmatrix}.$
\par
The Lie algebra ${\mathfrak gl}(n+2,\mathbb{R})$ is (as the vector
space) the direct sum of these four subspaces.  The commutator
relations between the basis introduced above are
\begin{multline}
[E^{i}_{j},E^{l}_{k}]=\delta^{l}_{j}E^{i}_{k}-\delta^{i}_{k}E^{l}_{j},\
[E^{i}_{j},1^{k}]=\delta^{k}_{j}1^{i},\
[E^{i}_{j},1_{k}]=-\delta_{k}^{i}1_{j},\\ [E_{**},1^{k}]=-1^{k},\
[E_{**},1_{k}]=1_{k},\
[1^{i},1_{j}]=E^{i}_{j}-\delta^{i}_{j}E_{**}.
\end{multline}

In the complexification ${\mathfrak gl}(n+2,\mathbb{C})$ of this
Lie algebra consider the {\bf real} Lie subalgebra ${\hat
{gl}}(n+1,\mathbb{R})$ of matrices \begin{equation}
\begin{pmatrix}
A & i{\bar f}\\
i{\bar e} & t
\end{pmatrix}
\end{equation}
where $A\in {\mathfrak gl}(n+1,\mathbb{R});\ {\bar f},{\bar e}\in
R^n ,\ t\in \mathbb{R}$.  The center of this Lie algebra is
one-dimensional, formed by the matrices proportional to the unit
matrix $I_{n+2}$. Factorization by this center defined the
epimorphism $\pi : {\hat {gl}}(n+1,\mathbb{R})\rightarrow {\hat
{sl}}(n+1,R)$ onto the Lie subalgebra of traceless matrices of the
form (23.33).\par

Define the embedding $j:{\mathfrak iso}_{\tG}\rightarrow
{\mathfrak gl}(n+2,\mathbb{C})$ as follows
\[
\sqrt{2}X_{k}\rightarrow i1_{k},\ \sqrt{2}D^{l}\rightarrow
i1^{l},\ Q^{i}_{j}\rightarrow E^{i}_{j}.
\]

It is easy to see that image of ${\mathfrak iso}_{\tG}$ under this
embedding is the subspace of ${\hat {\mathfrak
gl}}(n+1,\mathbb{R}).$ The commutator relations between the basic
elements $X_k,D^s,Q^{i}_{j}$ are preserved except for the last
one.  We have
\begin{multline}
j([\sqrt{2}D^{l},\sqrt{2}X_{k}])=j(-Q^{l}_{k}-\delta^{l}_{k}Q
)=-E^{l}_{k}-I_{n+1}=-\left(
E^{l}_{k}-\delta^{l}_{k}E_{**}+I_{n+2}\right)=\\
[i1^{l},i1_{k}]-I_{n+2}=
[j(\sqrt{2}D^{l}),j(\sqrt{2}X_{k})]-I_{n+2}.
\end{multline}
Combining the embedding $j$ with the projection $\pi : {\hat
{\mathfrak gl}}(n+1,\mathbb{R})\rightarrow {\hat {\mathfrak
sl}}(n+1,\mathbb{R})$ we get {\bf the isomorphism of Lie algebras}
\[
{\mathfrak iso}_{\tG}\simeq {\hat {\mathfrak sl}}(n+1,\mathbb{R}).
\]
The Lie algebras ${\mathfrak sl}(n+1,\mathbb{R})$ and ${\hat
{\mathfrak sl}}(n+1,\mathbb{R})$ are dual to one another in the
complex Lie algebra ${\mathfrak sl}(n+1,\mathbb{C})$ relative to
the involution
\[
\begin{pmatrix}
A & {\bar f}\\
{\bar e} & t
\end{pmatrix} \rightarrow \begin{pmatrix}
A & -{\bar f}\\
-{\bar e} & t
\end{pmatrix}
\]
and isomorphic as Lie algebras (see \cite{WG}).  Combining all the
calculations above we get the proof of the following

\begin{theorem}
The Lie algebra ${\mathfrak iso}_{\tG}\simeq {\mathfrak
sl}(n+1,\mathbb{R})$ of Killing vector fields of metric $\tG$ is
(as the vector space) the linear sum
\[
{\mathfrak iso}_{\tG}= {\mathfrak q}\oplus {\mathfrak d}\oplus
{\mathfrak x }
\]
of subalgebras
\[
1)\ {\mathfrak
q}=<Q^{i}_{j}=x^{i}\partial_{x^{j}}-p_{j}\partial_{p_{i}}>,
\]
with the commutator relations
\[
[Q^{i}_{j},Q^{p}_{k}]=\delta^{p}_{j}Q^{i}_{k}-\delta^{i}_{k}Q^{p}_{j},
\]
Subalgebra $\mathfrak q$ is isomorphic thus, to ${\mathfrak
gl}(n+1,\mathbb{R})$,\par

the abelian subalgebra
\[
2)\ {\mathfrak d}=<D^{i}=
\frac{x^{i}}{2}Q+(1-\frac{1}{2}(x^{l}p_{l}))\partial_{p_{i}},>
\]
where
$Q=\sum_{i}Q^{i}_{i}=\sum_{i}(x^{i}\partial_{x^{i}}-p_{i}\partial_{p_{i}})$
is the generator of hyperbolic rotation $H_{t}:(p,x)\rightarrow
(e^t p, e^{-t} q)$.
\par
and abelian subalgebra
\[
3)\ {\mathfrak x }=<X_{s}=\frac{\partial}{\partial x^{s}}>,
\]
\par
Generators $Q^{i}_{j},X_{i},D^{j}$ satisfy to the following
commutator relations
\begin{equation}
[Q^{i}_{j},X_{s}]=-\delta^{i}_{s}X_{j};\
[Q^{i}_{j},D^{s}]=\delta^{s}_{j}D^{i};\
[X_{s},D^{i}]=\frac{1}{2}Q^{i}_{s} +\frac{1}{2}\delta^{i}_{s}Q.
\end{equation}
\par
Vector fields $Q^{i}_{j},X_{i},D^{j}$ are Hamiltonian with
Hamiltonian functions
\[
H_{Q^{i}_{j}}=-x^i p_j;\ H_{X_k}=-p_k;\ H_{D^s}=x^s
(1-\frac{<{\bar x},{\bar p}>}{2}).
\]
\end{theorem}

\end{document}